\documentclass[11pt]{article}
\usepackage{amssymb}
\usepackage{epsfig}
\usepackage{amsmath}

\def\F{\Phi}

\newtheorem{theorem}{Theorem}

\newtheorem{lemma}[theorem]{Lemma}

\newtheorem{conjecture}[theorem]{Conjecture}
\newtheorem{remark}[theorem]{Remark}

\newenvironment{proof}{\noindent{\bf Proof. }}{\hfill$\square$\medskip}

\newcommand{\gap}[1]{\hspace{#1in}}
\newcommand{\brac}[1]{\left(#1\right)}
\newcommand{\bfrac}[2]{\brac{\frac{#1}{#2}}}
\newcommand{\beq}[1]{\begin{equation}\label{#1}}
\newcommand{\eeq}{\end{equation}}
\newcommand{\blem}[1]{\begin{lemma}\label{#1}}
\newcommand{\elem}{\end{lemma}}
\newcommand{\bth}[1]{\begin{theorem}\label{#1}}
\newcommand{\enth}{\end{theorem}}
\newcommand{\brem}[1]{\begin{remark}\label{#1}}
\newcommand{\erem}{\end{remark}}
\def\hD{\hat{D}}
\def\cE{{\cal E}}
\def\cA{{\cal A}}

\def\cC{{\cal C}}

\def\a{\alpha}

\def\d{\delta}
\def\D{\Delta}
\def\e{\varepsilon}
\def\f{\phi}
\def\g{\gamma}

\def\k{\kappa}

\def\th{\theta}

\def\l{\lambda}

\def\p{\pi}

\def\r{\rho}

\def\t{\tau}
\def\om{\omega}

\def\Th{\Theta}

\def\whp{{\bf whp}}

\newcommand{\rdup}[1]{\left\lceil #1 \right\rceil }

\newcommand{\set}[1]{\left\{#1\right\}}

\newcommand{\proofend}{\hspace*{\fill}\mbox{$\Box$}}

\def\E{{\sf E}}
\def\Pr{{\sf P}}

\newcommand{\ignore}[1]{}
\def\hp{\hat{p}}
\def\O1{\Omega(1)}
\newcommand{\card}[1]{\left|#1\right|}
\newcommand{\ts}[1]{#1^{(t)}}
\newcommand{\tsm}[1]{#1^{(t-1)}}
\newcommand{\tz}[1]{#1^{(0)}}
\newcommand{\tsd}[1]{#1^{(t')}}
\newcommand{\tss}[1]{#1^{(t+1)}}

\begin{document}
\author{Alan Frieze\thanks{Department of Mathematical Sciences, Carnegie Mellon University, Pittsburgh PA 15213. Supported in
part by NSF Grant DMS-0753472}\quad \quad \quad  Dhruv
Mubayi\thanks{Department of Mathematics, Statistics, and Computer
Science, University of Illinois, Chicago, IL 60607. Supported in
part by NSF Grant DMS-0653946}}

\title{Coloring simple hypergraphs}
\maketitle

\begin{abstract}
Fix an integer $k \ge 3$. A $k$-uniform hypergraph is simple if
every two edges share at most one vertex.  We prove that there is a
constant $c$ depending only on $k$ such that every simple
$k$-uniform hypergraph $H$ with maximum degree $\D$ has chromatic
number satisfying
$$\chi(H) <c \, \left(\frac{\D}{\log \D}\right)^{\frac{1}{k-1}}.$$
This implies a classical result of
Ajtai-Koml\'os-Pintz-Spencer-Szemer\'edi and its strengthening due
to Duke-Lefmann-R\"odl.  The result is sharp apart from the constant
$c$.

\end{abstract}
\section{Introduction}
Hypergraph coloring has been studied for almost 50 years, since
Erd\H os' seminal results on the minimum number of edges in uniform
hypergraphs that are not 2-colorable. Some of the major tools in
combinatorics have been developed to solve problems in this area,
for example, the local lemma and the nibble or semirandom method.
Consequently, the subject enjoys a prominent place among basic
combinatorial questions.

Closely related to coloring problems are questions about the
independence number of hypergraphs. An easy extension  of Tur\'an's
graph theorem shows that a $k$-uniform hypergraph with $n$ vertices
and average degree $d$ has an independent set of size at least $c
n/d^{1/(k-1)}$, where $c$ depends only on $k$.  If we impose local
constraints on the hypergraph, then this bound can be improved. An
$i$-cycle in a $k$-uniform hypergraph is a collection of $i$
distinct edges spanned by at most $i(k-1)$ vertices.  Say that a
$k$-uniform hypergraph has girth at least $g$ if it contains no
$i$-cycles for $2 \le i <g$. Call a $k$-uniform hypergraph simple if
it has girth at least 3. In other words, every two edges have at
most one vertex in common. Throughout this paper we will assume that
$k\ge 3$ is a fixed positive integer.

Ajtai-Koml\'os-Pintz-Spencer-Szemer\'edi \cite{AKPSS} proved the
following fundamental result that strengthened the bound obtained by
Tur\'an's theorem above.

\begin{theorem} {\bf (\cite{AKPSS})} \label{ind0}
 Let $H=(V,E)$ be a  $k$-uniform hypergraph of girth at least 5 with
maximum degree $\D$. Then it has an independent set of size at least
$$c\, n \left(\frac{\log \D}{\D}\right)^{1/(k-1)}$$
where $c$ depends only on $k$.
\end{theorem}
 Spencer conjectured that  Theorem \ref{ind0} holds even for simple
hypergraphs, and this was later proved by Duke-Lefmann-R\"odl
\cite{DLR}. Theorem \ref{ind0} has proved to be a seminal result in
combinatorics, with many applications. Indeed, the result was first
proved for $k=3$ by Koml\'os-Pintz-Szemer\'edi \cite{KPS} to
disprove the famous Heilbronn conjecture, that among every set of
$n$ points in the unit square, there are three points that form a
triangle whose area  is at most $O(1/n^2)$. For applications of
Theorem \ref{ind0} to coding theory or combinatorics,   see \cite{L}
or \cite{KMRT}, respectively.

The goal of this paper is to prove a result that is stronger than
Theorem \ref{ind0} (and also the accompanying result of \cite{DLR}).
Since the proof of our result does not use Theorem \ref{ind0}, it
gives an alternative proof of all the applications of Theorem
\ref{ind0} as well. Our main result states that not only can one
find an independent set of the size guaranteed by Theorem
\ref{ind0}, but in fact that the entire vertex set can be
partitioned into independent sets with this average  size.
 Recall that the chromatic number $\chi(H)$ of $H$ is
the minimum number of colors needed to partition the vertex set so
that no edge is monochromatic.

\begin{theorem}\label{th1a}
Fix $k \ge 3$. Let $H=(V,E)$ be a simple $k$-uniform hypergraph with
maximum degree $\D$. Then
$$\chi(H)<c{\bfrac{\D}{\log\D}} ^{\frac{1}{k-1}}$$
where $c$ depends only on $k$.
\end{theorem}

It is shown in \cite{BFM} that Theorem \ref{th1a} is sharp apart
from the constant $c$.
 In order to prove Theorem \ref{th1a} we will first
prove the following slightly weaker result. A triangle in a
$k$-uniform hypergraph is a 3-cycle that contains no 2-cycle.   In
other words, it is a collection of three sets $A, B, C$ such that
every two of these sets have nonempty intersection, and $A \cap B
\cap C=\emptyset$.

\begin{theorem}\label{th0}
Fix $k \ge 3$.  Let $H=(V,E)$ be a simple triangle-free $k$-uniform
hypergraph with maximum degree $\D$. Then
$$\chi(H)<c{\bfrac{\D}{\log\D}}^{\frac{1}{k-1}},$$
where $c$ depends only on $k$.
\end{theorem}

The proof of Theorem \ref{th0} rests on several major developments
in probabilistic combinatorics over the past 25 years. Our approach
is inspired by Johansson's breakthrough result on graph coloring,
which proved Theorem \ref{th0} for $k=2$.

 The proof technique, which has been termed
the semi-random, or nibble method, was first used by R\"odl
(inspired by earlier work in \cite{AKPSS, KPS}) to confirm the Erd\H
os-Hanani conjecture about the existence of asymptotically optimal
designs. Subsequently, Kim~\cite{Kim} (see also Kahn \cite{Kah})
proved Johansson's theorem for graphs with girth five and then
Johansson proved his result.  The approach used by Johansson for the
graph case is to iteratively color a small portion of the (currently
uncolored) vertices of the graph, record the fact that a color
already used at $v$ cannot be used in future on the uncolored
neighbors of $v$, and continue this process until the graph induced
by the uncolored vertices has small maximum degree. Once this has
been achieved, the remaining uncolored vertices are colored using a
new set of colors by the greedy algorithm.

In principle our method is the same, but there are several
difficulties we  encounter.  The first, and most important, is that
our coloring algorithm is not  as simple. A proper coloring of a
$k$-uniform hypergraph allows as many as $k-1$ vertices of an edge
to have the same color, indeed, to obtain optimal results one must
permit this. To facilitate this, we introduce a collection of $k-1$
different  hypergraphs at each stage of the algorithm whose edges
keep track of coloring restrictions. Keeping track of these
hypergraphs requires controlling more parameters during the
iteration and dealing with some more lack of independence and this
makes the proof more complicated.

In an earlier paper \cite{FM}, we had carried out this program for
$k=3$.  Several technical ideas incorporated in the current proof
can be found there. Because the notation in \cite{FM} is slightly
simpler than that in the current paper, the reader interested in the
technical details of our proof may want to familiarize him or
herself with \cite{FM} first (although this paper is entirely self
contained).

The implication Theorem \ref{th0} $\rightarrow$ Theorem \ref{th1a}
forms a much shorter (but still nontrivial) part of this paper (See
Section 2). Our proof uses a recent concentration result of Kim and
Vu \cite{KV} together with some additional ideas similar to those
from Alon-Krivelevich-Sudakov \cite{AKSud}. The approach here is to
partition the vertex set of a given simple hypergraph into some
number of parts, where the hypergraph induced by each of the parts
is triangle-free. Once this has been achieved, each of the parts is
colored using Theorem \ref{th0}.

Finally, we remark that our proof of Theorem \ref{th0} also gives
the same upper bound for list chromatic number, although we phrase
it only for chromatic number. On the other hand, we are not able to
prove Theorem \ref{th1a} for list chromatic number, since list
chromatic number is not additive in the sense described in the
previous paragraph.  We end with a conjecture posed in \cite{FM},
which states that we may replace the hypothesis "simple" with
something much weaker.

\begin{conjecture} {\bf (\cite{FM})} \label{conj}
Let $F$ be a  $k$-graph.  There is a constant $c_F$ depending only
on $F$ such that  every $F$-free $k$-graph with maximum degree
$\Delta$ has chromatic number at most $c_F(\Delta/\log
\Delta)^{1/{(k-1)}}$.
 \end{conjecture}
Conjecture \ref{conj} appears to be out of reach using current
methods. For example, the special case $k=2$ and $F=K_4$ remains
open and would imply an old conjecture of \cite{AEKS}.

Throughout this paper,  we will assume that $\D$ is sufficiently
large that all implied inequalities hold true.  Any asymptotic
notation is meant to be taken as $\D \to \infty$.
\section{Simple hypergraphs}
In this section we will prove that Theorem \ref{th0} $\rightarrow$
Theorem \ref{th1a}

Let $H=(V,E)$ be a simple $k$-uniform hypergraph. For $v\in V$ let
its neighbor set $N_H(v)$ be defined by $N_H(v)=\set{x:\;\exists
S\in E\ s.t.\ \{v,x\} \subset S}$. Let $d_H(v)$ denote the degree of
$v$ so that $|N_H(v)|=(k-1)d_H(v)$.


A pair $x,y\in N_H(v)$ is said to be {\em covered} if there exists
$S \in E$ that contains both $x$ and $y$ but not $v$. Note that $H$
simple implies that in this case no edge contains all of $v,x,y$.

Recall that $k$ is a fixed. Let $\e=\e(k)$ be a sufficiently small
positive constant depending only on $k$. Theorem \ref{th1a} will
follow from Theorem \ref{th0} and the following two lemmas:
\begin{lemma}\label{lem1}
Fix $k \ge 3$. Let $H=(V,E)$ be a simple $k$-uniform hypergraph with
maximum degree $\D$. Let
$m=\left\lceil\D^{\frac{2}{3k-4}-\e}\right\rceil$. Then there exists
a partition of $V$ into subsets $V_1,V_2,\ldots,V_m$ such each
induced subhypergraph $H_i=H[V_i],\,i=1,2,\ldots,m$ has the
following properties:
\begin{description}
 \item[(a)] The maximum degree $\D_i$ of $H_i$ satisfies $\D_i\leq 2\D/m^{k-1}$.
 \item[(b)] If $v\in V_i$ then its $H_i$-neighborhood $N_i(v)$ contains at most $k^2\D^2/m^{3k-4}$ covered pairs.
(Here we mean covered w.r.t. $H_i$).
\end{description}
\end{lemma}

\begin{lemma}\label{lem2}
Fix $k \ge 3$ and let $\d$ be a sufficiently small positive constant
depending
 on $k$.  Let $L=(V,E)$ be a simple $k$-uniform hypergraph with
maximum degree at most $d$. Suppose that each vertex neighborhood
$N_L(v)$ contains at most $d^\d$ covered pairs. Let
$\ell=d^{\frac{1}{k-1}-\d}$. Then there exists a partition of $V$
into subsets $W_1,W_2,\ldots,W_{\ell_1},\,\ell_1=O(\ell)$ such that
for each $1\leq j\leq \ell_1$, the hypergraph $L_j=L[W_j]$ has the
following properties:
\begin{description}
 \item[(a)] The maximum degree $d_j$ of $L_j$ satisfies $d_j\leq 2d/\ell^{k-1}$.
 \item[(b)] $L_j$ is triangle-free.
\end{description}
\end{lemma}
\proofend

\subsection{Proof of Theorem \ref{th1a}}
Our proof can be thought of as a nibble argument, involving two
iterations, given by Lemmas \ref{lem1} and \ref{lem2}. Suppose that
$H$ is a simple $k$-uniform hypergraph with maximum degree $\D$.
Apply Lemma \ref{lem1} to obtain $H_1, \ldots, H_m$ that satisfy the
conclusion of the lemma. Now fix $1 \le i \le m$ and let $L=H_i$.
Lemma \ref{lem1} part (a) implies that $\D(L)\le 2\D/m^{k-1}$. Hence
we may apply Lemma \ref{lem2} to $L$ with $d=2\D/m^{k-1}$ and
$\d=\e(3k-4)^2/(k-2)$. By Lemma \ref{lem1} part (b), each
neighborhood $N_L(v)$ contains at most $k^2\D^2/m^{3k-4}$ covered
pairs.  Now
$$k^2\frac{\D^2}{m^{3k-4}}\le k^2\D^{\e(3k-4)}=k^2\D^{\frac{\d(k-2)}{3k-4}}<d^\d.$$
We may therefore apply Lemma \ref{lem2} with
$\ell=d^{1/{(k-1)}-\d}$.  Together with Theorem \ref{th0} we obtain
$$
\chi(L)\le \sum_{j=1}^{\ell}
\chi(L_j)<O\brac{\ell\bfrac{d/\ell^{k-1}}{\log(d/\ell^{k-1})}^{\frac{1}{k-1}}}=O\brac{
\bfrac{d}{\log({d^{(k-1)\d}})}^{\frac{1}{k-1}}}=
O\brac{\bfrac{d}{\log d}^{\frac{1}{k-1}}}. $$ Since this holds for
each $H_i$ we obtain,
$$\chi(H)\le \sum_{i=1}^m \chi(H_i)
<O\brac{m\bfrac{d}{\log d}^{\frac{1}{k-1}}}=
O\brac{m\bfrac{\D/m^{k-1}}{\log(\D/m^{k-1})}^{\frac{1}{k-1}}}=O\brac{\bfrac{\D}{\log\D}^{\frac{1}{k-1}}}.$$

\subsection{Kim-Vu concentration}
We will need the following very
useful concentration inequality \eqref{KimVu} from Kim and Vu
\cite{KV}: Let $\Upsilon=(W,F)$ be a hypergraph of rank $s$, meaning
that each $f\in F$ satisfies $|f|\leq s$. Let
$$Z=\sum_{f\in F}\prod_{i\in f}z_i$$
where the $z_i,i\in W$ are independent random variables taking
values in $[0,1]$. For $A\subseteq W,|A|\leq s$ let
$$Z_A=\sum_{\substack{f\in F\\f\supseteq A}}\prod_{i\in f\setminus A}z_i.$$
Let $M_A=\E(Z_A)$ and $M_j=\max_{A, |A|\geq j}M_A$ for $j\ge0$.
There exist positive constants $a=a_s$ and $b=b_s$ such that for any
$\l>0$,
\begin{equation}\label{KimVu}
\Pr(|Z-\E(Z)|\geq a\l^{s}\sqrt{M_0M_1})\leq b|W|^{s-1} e^{-\l}.
\end{equation}

\subsection{Proof of Lemma \ref{lem1}}

We will use the local lemma in the form below.

\begin{theorem} {\bf (Local Lemma)} \label{ll}
Let $\cA_1, \ldots, \cA_n$ be events in an arbitrary  probability
space.  Suppose that each event $\cA_i$ is mutually independent of a
set of all the other events $\cA_j$ but at most $d$, and that
$P(\cA_i)<p$ for all $1 \le i \le n$.  If $ep(d+1)<1$, then with
positive probability, none of the events $\cA_i$ holds.
\end{theorem}

We will partition $V$ randomly
into $m$ parts of size $\sim|V|/m$ and use the local lemma to show
the existence of a partition. We make the partition by assigning a
random number in $[m]$ to each $v\in V$.

Fix $v\in V$. To simplify notation, condition on $v\in V_1$. Let
$A_v$ be the event that (a) fails at $v$ i.e. that $v$ has degree
greater than $2\D/m^{k-1}$ in the hypergraph $H_1$.

Let $B_v$ be the event that its neighborhood in $H_1$ contains more
than $k^2\D^2/m^{3k-4}$ covered pairs.

Each of these events is mutually independent of a set of all other
events but at most $O(\D^4)$. We will show that\\
$\Pr(A_v),\Pr(B_v) =O(\D^{-5})$ and this is clearly sufficient for
the application of the local lemma.

Let $d_v$ be the degree of $v$ in $H_{1}$. Then $d_v$ has a
distribution that is dominated by the binomial distribution
$Bin(\D,1/m^{k-1})$. It follows from the Chernoff bounds that
$$\Pr\brac{d_v\ge 2\D/m^{k-1}}\leq e^{-\D/(3m^{k-1})}=e^{-\D^{\frac{k-2}{3k-4}-(k-1)\e+o(1)}}\leq \D^{-5}$$
and this disposes of $A_v$.

Our goal now is  to bound $\Pr(B_v)$.
 For a vertex $x \in N_H(v)$, let
$T_v(x)$ denote the unique $(k-1)$-set containing $x$ such that
$T_v(x) \cup \{v\} \in E$.  A covered pair $x,y\in N_H(v)$ will
remain as a covered pair in $N_{H_1}(v)$ iff $S_{x,y}=T_v(x) \cup
T_v(y)  \cup T\subseteq V_1$ where $T$ is the unique $(k-2)$-set
such that $T \cup \{x,y\} \in E$. Let $S_1,S_2,\ldots,S_r,\,r\leq
(k-1)^2\binom{\D}{2}$ be an enumeration of the $(3k-4)$-tuples
$S_{x,y}$ as $\set{x,y}$ ranges over the covered pairs in $N_H(v)$.

We will use the concentration inequality (\ref{KimVu}). The edges of
our hypergraph $(W,F)$ are $S_1,S_2,\ldots,S_r$ and if $x\in W$ then
$z_x$ is an independent $\set{0,1}$ Bernoulli random variable with
$\Pr(z_x=1)=1/m$. Note that $|W|\leq k\D^2$.

Let $Z_v$ denote the number of covered pairs in $N_{H_1}(v)$. There
is a 1-1 correspondence between covered pairs and the $S_i$.
Therefore \beq{Zv} \mu = \E(Z_v) = \frac{r}{m^{3k-4}} \leq
\frac{(k-1)^2\D^2}{2m^{3k-4}}. \eeq We now have to estimate $M_1$.

For each set $A \subset W$, let $Y_A$ denote the number of edges of
$F$ containing $A$.

{\bf Claim.} $|Y_A| =O(\D)$ if $|A| \le k-1$ and $|Y_A|=O(1)$ if
$|A|\geq k$.

\begin{proof} Suppose that $A \subset S \in F$. Then $S$ can be written as
$T_v(x) \cup T_v(y) \cup T$, for some $x,y \in N_H(v)$. We will
count the number of $S$ containing $A$ by the number of $T_v(x)$'s
and $T$'s. First, the number of $S$ where both $T_v(x)$ and $T_v(y)$
have a vertex in $A$ is at most $5k^4$ by the following argument.
There are $\binom{|A|}{2}<5k^2$ choices for the two intersection
points, these points uniquely determine $T_v(x)$ and $T_v(y)$, and
there are at most $|T_v(y)||T_v(x)|=(k-1)^2$ possible covered pairs,
each of which determines $T$ uniquely (if $T$ exists).

Now suppose that $A \cap T_v(y)=\emptyset$. If $|A| \geq
k$, then $A$ must contain a vertex from $T$ and a
(different) vertex from $T_v(x)$.  There are at most $9k^2$ choices
for these  two vertices.  For each of these choices, $T_v(x)$ is
determined uniquely.  For each vertex in $T_v(x)$ and the chosen
vertex of $T \cap A$, there is at most one choice for $T$ (since $H$
is simple), hence the number of choices for $T$ is at most $k-1$.
Having chosen $T$, there are at most $k-1$ choices for $T_v(y)$.
Altogether, there are at most $9k^4$ choices for $S$. We conclude
that if $|A|\geq k$, then
$$|Y_A| \le 5k^4+9k^4=O(1).$$
If $|A| \le k-1$, the argument above still applies unless either $A
\subset T_v(x)$ or $A \subset T$. In either case, there are at most
$k\D$ ways of choosing the other part of $S \setminus T_v(y)$ and at
most $k$ ways of choosing $T_v(y)$.  Thus $|Y_A|=O(\D)$ as claimed.
\end{proof}

The probability of choosing each vertex in $S\setminus A$ is $1/m$,
so for given $A$, the probability of a particular $S \supset A$ is
$(1/m)^{3k-4-|A|}$. The Claim now implies that for $1 \le |A| <
3k-4$,
$$M_A \le \max\left\{O\left(\frac{\D}{m^{2k-3}}\right), O
\left(\frac{1}{m}\right)\right\}.$$
By our choice of $m$, it follows
that if $\e$ is sufficiently small then
$$M_1= O(\D^{-1/(3k)}).$$

 The choice of $m$ also gives
$$M_0=\max\{\mu, M_1\} \le \frac{k^2\D^2}{m^{3k-4}}.$$
It follows that if we take $a\lambda^{3k-4}= \frac{k^2\D^2}{m^{3k-4}(M_0M_1)^{1/2}}$ then  $\lambda>\D^{\delta_k}$ where $\delta_k>0$.
Now \eqref{KimVu} implies that
$$\Pr(B_v)\leq \Pr(Z_v-\E(Z_v)\geq k^2 \D^2/m^{3k-4})\leq b(k\D^2)^{3k}e^{-\lambda}\leq \D^{-5}.$$
This completes the proof of Lemma \ref{lem1}. \proofend
\subsection{Proof of Lemma \ref{lem2}}
This part follows an approach taken in Alon, Krivelevich and Sudakov
\cite{AKSud}. We will first partition $V$ randomly into $\ell$ parts
$V_1,V_2,\ldots,V_{\ell}$ of size $\sim|V|/\ell$ and use the local
lemma  to show the existence of a partition satisfying certain
properties. To simplify notation, condition on $v\in V_1$.

For $v\in V$ let $A_v$ be the event that (a) fails at $v$ i.e. that
$v$ has degree greater than $2d/\ell^{k-1}$ in $L_1$.

Let $B_v$ be the event that $N_{L_1}(v)$ contains at least $200k^2$
covered pairs w.r.t. $L_1$.

Each of these events is mutually independent of a set of all other
events but at most $O(d^4)$.  We will show that\\
$\Pr(A_v),\Pr(B_v) =O(d^{-5})$ and this is clearly sufficient for
the application of the local lemma.

Let $d_v$ be the degree of $v$ in $L_{1}$. Then $d_v$ has a
distribution that is dominated by the binomial distribution
$Bin(d,1/\ell^{k-1})$. It follows from the Chernoff bounds that
$$\Pr\brac{d_v\ge 2d/\ell^{k-1}}\leq e^{-d/(3\ell^{k-1})}=e^{-d^{(k-1)\d}/3}\leq d^{-5}$$
and this disposes of the $A_v$.

If $B_v$ fails then either
\begin{description}
\item[(i)] There exists a vertex $w\in N_L(v)$ such that $w$ is in at least $10k$ covered pairs of
$N_{L_1}(v)$, or
\item[(ii)] $N_{L_1}(v)$ contains at least $10k$ pair-wise disjoint covered pairs.
\end{description}
Now
\begin{eqnarray*}
\Pr((i))&\leq&kd\binom{d^\d}{10k}\ell^{-10k}\leq d^{-5}\\
\Pr((ii))&\leq&\binom{d^{2\d}}{10k}\ell^{-20k}\leq d^{-5}
\end{eqnarray*}
and this  disposes of the $B_v$.

So, assume that none of the events $A_v,B_v$ occur. We show now that
we can partition each $V_j$ into at most $400k^2+1$ sets, each of
which induces a triangle free hypergraph. Consider the digraph $D_1$
with vertex set $V_1$ and an edge directed from $v\in V_1$ to each
vertex of each of the at most $200k^2$ covered pairs in
$N_{L_1}(v)$. $D_1$ has maximum out-degree $400k^2$ and so its
underlying graph $G_1$ is $400k^2$-degenerate and so it can be
properly colored with $400k^2+1$ colors. Partition $V_1$ into color
classes $W_1,W_2,\ldots,W_{400k^2+1}$. We claim that for each $s$,
the hypergraph $L[W_s]$ induced by $W_s$ is triangle-free. Suppose
then that there is a triangle $v \cup T_v(x), v \cup T_v(y), T$
inside $L[W_s]$, where $T$ contains both $x$ and $y$. Then
$\set{x,y}$ is a covered pair for $v$ and by construction $v$ and
$x$ are not in the same $W_s$, contradiction. \proofend

\section{Triangle-free hypergraphs}
In this section, which forms the bulk of the paper, we will prove
Theorem \ref{th0}.

\subsection{Local Lemma}

The driving force of our upper bound argument, both in the
semi-random phase and the final phase, is the Local Lemma. Note that
the Local Lemma immediately implies that every $k$-graph with
maximum degree $\Delta$ can be properly colored with at most
$\rdup{4\Delta^{1/(k-1)}}$ colors. Indeed, if we color each vertex
randomly and independently with one of these colors, the probability
of the event $\cA_e$, that an edge $e$ is monochromatic, is at most
$\frac{1}{4^{k-1}\Delta}$.  Moreover $\cA_e$ is independent of all
other events $\cA_f$ unless $|f\cap e|>0$, and the number of $f$
satisfying this is less than $k\Delta$. We conclude that there is a
proper coloring.

\subsection{Coloring Procedure}
In the rest of the paper, we will
prove the upper bound in Theorem \ref{th0}.  Suppose that $k\ge 3$
is fixed and $H$ is a simple triangle-free $k$-graph with maximum
degree $\Delta$.

 Let $V$ be the vertex set of $H$. As usual, we write $\chi(H)$ for
the chromatic number of $H$.  Let $\e$ be a sufficiently small fixed
positive constant (depending only on $k$).  Let
$$\om=\frac{\e^2\log\D}{100\times k^{2k+1}}$$
and set
$$q=\left\lceil\frac{\D^{1/(k-1)}}{\om^{1/(k-1)}}\right\rceil.$$
Note that $q <c(\Delta/\log \Delta)^{1/(k-1)}$ where $c$ depends
only on $k$.

We color $V$ with $2q$ colors and therefore show that
$$\chi(H)\leq 2c\bfrac{\Delta}{\log \Delta}^{1/(k-1)}.$$
We use the first $q$ colors to color $H$ in rounds and then
use the second $q$ colors to color any vertices not colored by this process.

Our algorithm for coloring in rounds is semi-random.
At the beginning of a round certain parameters will satisfy certain properties, {\bf \eqref{A} -- \eqref{Dc}} below.
We describe a set of random choices for the parameters in the next round and
we use the local lemma to prove that there is a set of choices that preserves the required properties.

\begin{itemize}
\item $C=[q]$ denotes the set of available colors for the semi-random phase.
\item  $\ts{U}$: The set of vertices which are currently uncolored. ($\tz{U}=V$).
\item $\ts{H}$: The sub-hypergraph of $H$ induced by $\ts{U}$.
\item $\ts{W}=V\setminus \ts{U}$: The set of vertices that have been colored. We use the notation $\k$ to denote the color
of an item e.g. $\k(w),\,w\in \ts{W}$ denotes the color permanently assigned to $w$.
\item $\ts{H_i},\,2\leq i\leq k-1$: An edge-colored $i$-graph with vertex set $\ts{U}$. There is an edge
$u_1u_2\cdots u_i\in \ts{H_i}$ iff there are vertices
$u_{i+1},u_{i+2},\ldots,u_k\in \ts{W}$ and an edge $u_1u_2\cdots
u_k\in H$ with $\k(u_{i+1})=\k(u_{i+2})=\cdots=\k(u_k)$. For a fixed
$u_1u_2\cdots u_i$, this color is well defined because of the fact
that $H$ is simple. The edge $u_1u_2\cdots u_i$ is given the color
$\k(u_{i+1})$. (These hypergraphs are used to keep track of coloring
restrictions).
\item  $\ts{p_u}\in [0,1]^C$ for $u\in \ts{U}$: This is a vector of coloring probabilities. The $c$th coordinate
is denoted by $\ts{p}_u(c)$ and
$\tz{p_u}=(q^{-1},q^{-1},\ldots,q^{-1})$.
\end{itemize}
We can now describe the ``algorithm'' for computing
$\tss{U},H_i^{(t+1)}, \tss{p}_u$, given $\ts{U}, H_i^{(t)},
\ts{p}_u,$ for $u\in \ts{U}$: Let
$$\th=\frac{\e}{{\om}}= \frac{100\times k^{2k+1}}{\e\log\D}$$
where we recall that $\e$ is a sufficiently small positive constant.

For each $u\in \ts{U}$ and $c\in C$ we {\em tentatively activate}
$c$ at $u$ with probability $\th \ts{p}_u(c)$. A color $c$ is lost
at $u\in \ts{U}$, $\tss{p}_u(c)= 0$ and $\tsd{p}_u(c)=0$ for $t'>t$
if either

(i) there is an edge $uu_2\cdots u_k\in \ts{H}$ such that $c$ is
tentatively activated at $u_2,u_3,\ldots,u_k$ or

(ii) there is a $2\leq i\leq k-1$ and an edge $e=uu_2\cdots u_i\in
\ts{H}_i$ such that $c=\k(e)$ and $c$ is tentatively activated at
$u_2,u_3,\ldots,u_i$.

The vertex $u\in \ts{U}$ is given a permanent color if there is a
color tentatively activated at $u$ which is not lost due to the
above reasons. If there is a choice, it is made arbitrarily. Then
$u$ is placed into $\tss{W}$.

We fix
$$\hp=\frac{1}{\D^{1/(k-1)-\e}}.$$

We keep
$$\ts{p}_u(c)\leq \hp$$
for all $t,u,c$.

We let
$$\ts{B}(u)=\set{c:\;\ts{p}_u(c)=\hp}\qquad for\ all\ u\in V.$$
A color in $\ts{B}(u)$ cannot be used at $u$. The role of $\ts{B}(u)$ is clarified later.

Suppose that color $c$ has not been lost at $u$ prior to round $t$.
Let us compute the probability that $c$ is not lost at $u$ in round
$t$.  Since for each $u \in U^{(t)}$ and $c \in C$, the tentative
activation of $c$ at $u$ is done independently of all other
tentative activations, the probability that $c$ is not lost at $u$
due to (i) is
$$\prod_{uu_2\cdots u_k\in
\ts{H}}\brac{1-\prod_{j=2}^k \th\ts{p}_{u_j}(c)}=\prod_{uu_2\cdots
u_k\in \ts{H}}\brac{1-\th^{k-1}\prod_{j=2}^k\ts{p}_{u_j}(c)}.$$
Similarly the probability that $c$ is not lost at $u$ due to (ii) is
$$\prod_{i=2}^{k-1}\prod_{\substack{e=uu_2\cdots u_i\in \ts{H}_i\\ \k(e)=c}}\brac{1-\th^{i-1}\prod_{j=2}^i \ts{p}_{u_j}(c)}.$$
Consequently, the probability that $c$ is not lost at $u$ in round
$t$, given that it wasn't lost in any prior round is
\beq{2a}
\ts{q}_u(c)=\prod_{uu_2\cdots u_k\in \ts{H}}\brac{1-\th^{k-1}\prod_{j=2}^k\ts{p}_{u_j}(c)}
\prod_{i=2}^{k-1}\prod_{\substack{e=uu_2\cdots u_i\in \ts{H}_i\\
\k(e)=c}}\brac{1-\th^{i-1}\prod_{j=2}^i \ts{p}_{u_j}(c)}
\eeq
The
parameter $\ts{q}_u(c)$ is of great importance in our proof.

{\bf Coloring Procedure: Round $t$}
\begin{description}
\item[Make tentative random color choices:]\ \\
Independently, for all $u\in \ts{U}$, $c\in C$, let
\beq{gammadef}
\ts{\g}_u(c)=\begin{cases}
           1&Probability\ =\ \th \ts{p}_u(c)\\0&Probability\ =\ 1-\th \ts{p}_u(c)
           \end{cases}
\eeq
$$\ts{\Th}(u)=\set{c:\;\ts{\g_u}(c)=1}\ =\ \mbox{the set of colors tentatively activated at }u.$$
\item[Deal with color clashes:]
\begin{multline*}
\ts{L}(u)=\set{c:\;\exists uu_2\cdots u_k\in \ts{H},\, such\ that\ c\in \bigcap_{j=2}^k\ts{\Th}(u_j)}\cup\\
\set{c:\;\exists 2\leq i\leq k-1\ and\ e=uu_2\cdots u_i\in \ts{H}_i\ such\ that\ \k(e)=c\in  \bigcap_{j=2}^i\ts{\Th}(u_j)}
\end{multline*}
is the set of colors {\em lost} at $u$ in this round.
$$\ts{A}(u)=\tsm{A}(u)\cup\ts{L}(u).$$
\item[Assign some permanent colors:]\ \\ Let
$$\ts{\Psi}(u)=\ts{\Th}(u)\setminus (\ts{A}(u)\cup \ts{B}(u))
\mbox{ = set of activated colors that can be used at }u.$$

If $\ts{\Psi}(u)\neq \emptyset$ then choose $c\in \ts{\Psi}(u)$ arbitrarily. Let $\k(u)=c$.
\end{description}
{\bf Update parameters:}
\begin{description}
\item[(a)]
$$\tss{U}=\ts{U}\setminus \set{u:\;\ts{\Psi}(u)\neq \emptyset}.$$
\item[(b)]
$\tss{H}_i,\,2\leq i\leq k-1$ is the $i$-graph with vertex set
$\tss{U}$ and edge set
$$\{u_1u_2\cdots u_i:\;\exists u_{i+1},\ldots, u_k \in W^{(t+1)}
\hbox{ with } u_1\cdots u_k \in H \hbox{ and }
\k(u_{i+1})=\cdots=\k(u_k)=c\}.$$
 Edge $ u_1u_2\cdots
u_i$ has color $c$. ($H$ simple implies that this color is
well-defined).
\item[(c)] $\ts{p}_u(c)$ is replaced by $p_u^{(t+1)}(c)$ which is
either 0, $p_u^{(t)}(c)/q_u^{(t)}(c)$, or $\hp$  (note that the last
two are at least $p_u^{(t)}(c)$). Furthermore, if $u\in
\ts{U}\setminus \tss{U}$ then by convention $\tsd{p}_u=\tss{p}_u$
for all $t'> t$.

In order to decide which of these three values is taken by
$p_u^{(t+1)}(c)$, we perform a random experiment, where we replace
$p_u^{(t)}(c)$ by  a random value $p_u'(c)$. Based on the outcome of
this random experiment, we will decide on the value of
$p_u^{(t+1)}(c)$.  One of the key properties is
 \beq{puc1}
\E(p'_u(c))=\ts{p}_u(c).
\eeq
\vspace{.2in}
The update rule is as
follows: If $c\in \tsm{A}(u)$ then $\ts{p}_u(c)$ remains unchanged
at zero. Otherwise, let $\ts{\eta}_u(c)$ be a random variable with
$$\ts{\eta}_u(c)\in \set{0,1}\ and\ \Pr(\ts{\eta}_u(c)=1)=\ts{p}_u(c)/\hp,\ independently\ of\ other\ variables.$$
Then
\beq{pdash}
p_u'(c)=\begin{cases}
            \begin{cases} 0&c\in \ts{L}(u)\\ \frac{\ts{p}_u(c)}{\ts{q}_u(c)}&c\notin \ts{L}(u)
\end{cases}&\frac{\ts{p}_u(c)}{\ts{q}_u(c)}< \hp\qquad\qquad\qquad {\bf Case\ A}\\ \\
            \ts{\eta}_u(c)\hp&\frac{\ts{p}_u(c)}{\ts{q}_u(c)}\geq \hp.\qquad \qquad\qquad{\bf Case\ B}
             \end{cases}
\eeq
\end{description}
There will be
$$t_0=\e^{-1}\log\D\log\log\D\ rounds.$$

Before getting into the main body of the proof, we check
\eqref{puc1}.  First  observe that $q_u^{(t)}(c)$ is the probability
that $c \not\in L^{(t)}(u)$ given that $c \not\in A^{(t-1)}(u)$.

If $\ts{p}_u(c)/\ts{q}_u(c)< \hp$ then
$$\E(p'_u(c))=\ts{q}_u(c)\frac{\ts{p}_u(c)}{\ts{q}_u(c)}=\ts{p}_u(c).$$
If $\ts{p}_u(c)/\ts{q}_u(c)\geq \hp$ then
$$\E(p'_u(c))=\hp\frac{\ts{p}_u(c)}{\hp}=\ts{p}_u(c).$$
Note that once a color enters $\ts{B}(u)$, it will be in
$B^{(t')}(u)$ for all $t'\geq t$. This is because we update $p_u(c)$
according to Case B and now $\Pr(\ts{\eta}_u(c)=1)=1$. We arrange
things this way, because we want to maintain \eqref{puc1}. Then
because $\ts{p}_u(c)$ cannot exceed $\hp$, it must actually remain
at $\hp$. This could cause some problems for us if neighbors of $u$
had been colored with $c$. There might be an edge $e=uu_2\cdots u_k$
where $u_2,\ldots,u_k$ are (tentatively) colored $c$. We don't want
to raise $\ts{p}_u(c)$ and to keep it monotone, we can't allow it to
drop to zero. This is why $\ts{B}(u)$ is excluded in the definition
of $\ts{\Psi}(u)$ i.e. we cannot color $u$ with $c\in \ts{B}(u)$.
\subsection{Correctness of the coloring}
Observe that if color $c$
enters $\ts{A}(x)$ at some time $t$ then $\k(x)\neq c$ since
$A^{(i)}(x)\subseteq A^{(i+1)}(x)$ for all $i$. Suppose that some
edge $u_1u_2\cdots u_k$ is improperly colored by the above
algorithm. Suppose that $u_1,u_2,\ldots,u_k$ get colored at times
$t_1\leq t_2\leq \cdots\leq t_k$ and that $\k(u_j)=c$ for
$j=1,2,\ldots,k$. If $t_1=t_2=\cdots=t_{k-1}=t$ then $c\in
\ts{L}(u_k)$ and so $\k(u_k)\neq c$. If there exists $1\leq i\leq
k-2$ such that $t_i<t=t_{i+1}=\cdots =t_{k-1}$  then
$u_{i+1}u_{i+2}\cdots u_k$ is an edge of $\ts{H}_{k-i}$ and
$\k(u_{i+1}u_{i+2}\cdots u_k)=c$ and so $c\in \ts{L}(u_k)$ and again
$\k(u_k)\neq c$.
\subsection{Parameters for the problem}
We will now
drop the superscript $(t)$, unless we feel it necessary. It will be
implicit i.e. $p_u(c)=\ts{p}_u(c)$ etcetera. Furthermore, we use a
$'$ to replace the superscript $(t+1)$ i.e. $p_u'(c)=\tss{p}_u(c)$
etcetera.  The following are the main parameters that we need in the
course of the proof:

In what follows $u_1=u$ and $2\leq i\leq k-1$:
\begin{eqnarray*}
\Xi_e&=&\sum_{c\in C}\prod_{j=1}^kp_{u_j}(c)\qquad\ for\ edge\ e=u_1u_2\cdots u_k\ of\ {\ts{H}}.\\
\F_{u,i}&=&\sum_{c\in C}\ \sum_{\substack{e=uu_2\cdots u_i\in H_i\\ \k(e)=c}}\ \prod_{j=1}^ip_{u_j}(c)\\
{h}_u&=&-\sum_{c\in C}p_u(c)\log p_u(c).\\
d_i(u,c)&=&|\set{e:\;u \in e\in H_i\ and\ \k(e)=c}|\\
d_i(u)&=&\sum_{c\in C}d_i(u,c) =\ degree\ of\ u\ in\ H_i\\
d_{\ts{H}}(u)&=&|\{e:\; u \in e\in \ts{H}\}|=\ degree\ of\ u\ in\ \ts{H}\\
d(u)&=&d_2(u)+d_3(u)+\cdots+d_{k-1}+d_{\ts{H}}(u)
\end{eqnarray*}

It will also be convenient to define the following auxiliary
parameters:
\begin{eqnarray*}
\Xi_e(c)&=&\prod_{j=1}^kp_{u_j}(c)\qquad\ for\ edge\ e=u_1u_2\cdots u_k\ of\ {\ts{H}}.\\
\Xi_u&=&\sum_{e=uu_2\cdots u_k\in \ts{H}}\Xi_{e}\\
\Xi_u(c)&=&\sum_{uu_2\cdots u_k\in \ts{H}}\prod_{j=2}^k p_{u_j}(c)\\
{\F}_{u,i}(c)&=&\sum_{\substack{e=uu_2\cdots u_i\in H_i\\ \k(e)=c}}\
\prod_{j=2}^i p_{u_j}(c)
\end{eqnarray*}
This gives
\begin{eqnarray}
\Xi_u&=&\sum_{c\in C}p_u(c)\Xi_u(c)\label{N1}\\
\F_{u,i}&=&\sum_{c\in C}p_u(c)\F_{u,i}(c).\label{N2}
\end{eqnarray}

\subsection{Invariants}
We define a set of properties such that if
they are satisfied at time $t$ then it is possible to extend our
partial coloring and maintain these properties at time $t+1$. These
properties are now listed. They are only claimed for $u\in U$ and
they are easily verified for $t=0$.
\begin{eqnarray}
\card{1-\sum_c p_u(c)}&\leq&t\D^{-\e}.\label{A}\\
\nonumber\\
\Xi_{e}&\leq&\Xi_{e}^{(0)}+\frac{t}{\D^{1+\e}}\label{B}\\
&=&\frac{{\om}}{\D}+\frac{t}{\D^{1+\e}}\,,\qquad\forall e\in \ts{H}.\nonumber\\
\nonumber\\
\F_{u,i}&\leq&k^{2k-2i}\om(1-\th/3k)^t,\qquad 2\leq i\leq k-1.\label{Ca}\\
\nonumber\\
h_u&\geq&h^{(0)}_u-k^{2k}\e \sum_{\t=0}^t(1-\th/3k)^\t.\label{C}\\
\nonumber\\
d(u)&\leq&\brac{1-\th/2k}^t\D\label{D}.\\
\nonumber\\
d_i(u,c)&\leq &(1+2k\th)^t\D\hp^{k-i},\qquad 2\leq i\leq k-1.\label{Dc}
\end{eqnarray}
Equation \eqref{D} shows that after $t_0$ rounds we find that the maximum degree in the hypergraph
induced by the uncolored vertices satisfies
\begin{eqnarray}
\D(H^{(t_0)})&\leq&\brac{1-\th/2k}^{t_0}\D\nonumber\\
&\leq&e^{-\th t_0/2k}\D\nonumber\\
&<&e^{-100k^{2k}\log\log\D/2\e^2}\D\nonumber\\
&=& \frac{\D}{(\log\D)^{100k^{2k}/2\e^2}}.\label{DHt0}
\end{eqnarray}
and then the local lemma will show that the remaining vertices can
be colored with a set of
$4(\D/(\log\D)^{100k^{2k}/2\e^2})^{1/(k-1)}+1<q$ new colors.

The above invariants allow us to prove the following bounds:
By repeatedly using $(1-a)(1-b)\geq 1-a-b$ for $a,b\geq 0$ we see that
\beq{crude1}
q_u(c)\geq1-\th^{k-1}\Xi_u(c)-\sum_{i=2}^{k-1}\th^{i-1} \F_{u,i}(c).
\eeq

\subsection{Dynamics}
To prove \eqref{A} -- \eqref{Dc} we show that
we can find updated parameters such that
\begin{eqnarray}
\card{\sum_cp_u'(c)-\sum_cp_u(c)}&\leq&\D^{-\e}.\label{A1}\\
\nonumber\\
\Xi_{e}'&\leq& \Xi_{e}+\D^{-1-\e}.\label{B1}\\
\nonumber\\
\F_{u,i}'-\F_{u,i}&\leq&\binom{k-1}{i-1}
\th^{k-i}\Xi_u+\sum_{l=i+1}^{k-1}\binom{l-1}{i-1}\th^{l-i}\F_{u,l}\\
\nonumber & & -\th(1-k^{2k}\e)\F_{u,i}+\D^{-\e}, \qquad\qquad \qquad \qquad \qquad 2 \le i \le
k-1.
\label{B1a}\\
\nonumber\\
h_u-h_u'&\leq&k^{2k}\e(1-\th/3k)^t.\label{C1}\\
\nonumber\\
d'(u)&\leq&(1-\th/k)d(u)+\D^{2/3}.\label{D1}\\
\nonumber\\
d_i'(u,c)&\leq&d_i(u,c)+2k\th(1+2k\th)^t\D\hp^{k-i}, \quad 2 \le i \le
k-1.\label{duc}
\end{eqnarray}
\subsection{\eqref{A1}--\eqref{duc} imply \eqref{A}--\eqref{Dc}}
First let us show that \eqref{A1}--\eqref{duc} are enough to
inductively prove that \eqref{A}--\eqref{D} hold throughout.

{\bf Property \eqref{A}:} Trivial.

{\bf Property \eqref{B}:} Trivial.

{\bf Property \eqref{Ca}:} Fix $u$ and note that \eqref{B} and
\eqref{D} imply
\beq{eu}
\Xi_u\leq
\brac{\frac{{\om}}{\D}+t\D^{-1-\e}}d(u)\leq
{\om}(1-\th/2k)^t+\D^{-\e/2}.
\eeq
Therefore,
$$\F_{u,i}'-\F_{u,i}\leq\binom{k-1}{i-1}\th^{k-i}{\om}(1-\th/2k)^t
+\sum_{l=i+1}^{k-1}\binom{l-1}{i-1}\th^{l-i}\F_{u,l}-\th(1-k^{2k}\e)\F_{u,i}+\D^{-\e/3}$$
from \eqref{B1a} and \eqref{eu}.
Thus,
\begin{eqnarray*}
\F_{u,k-1}'&\leq& (k-1)\th\om(1-\th/2k)^t+(1-\th(1-k^{2k}\e))k^2\om(1-\th/3k)^t+\D^{-\e/3}\\
&\leq&k^2(1-\th/3k)^{t+1}\om\brac{\frac{\th(1-\th/2k)^t}{k(1-\th/3k)^{t+1}}
+
\frac{1-\th(1-k^{2k}\e)}{1-\th/3k}}+\D^{-\e/3}\\
&\leq&k^2(1-\th/3k)^{t+1}\om.
\end{eqnarray*}

Now for $i\leq k-2$, using $\F_{u,l}\leq k^{2k-2l}\om(1-\th/3k)^t$,
\begin{align*}
 &\F_{u,i}'\leq \binom{k-1}{i-1}\th^{k-i}\om(1-\th/2k)^t
+\frac{k^{2k}}{\th^i}\om(1-\th/3k)^t\sum_{l=i+1}^{k-1}\binom{l-1}{i-1}
\bfrac{\th}{k^2}^l\\
&\gap{2}+(1-\th(1-k^{2k}\e))k^{2k-2i}\om(1-\th/3k)^t+\D^{-\e/3}.
\end{align*}
Factoring out the first term $\binom{i+1}{i-1}(\th/k^2)^{i+1}$ in
the sum above we are left with a sum that can be upper bounded by
$\sum_{i=0}^{\infty}r^i$ where $0<r<\th(k-1)/k^2$.  Since
$\th=O(1/\log \Delta)$, this geometric series is upper bounded by
$1+\th/k$.  Consequently, $\F_{u,i}'$ is upper bounded by
\begin{align*}
&\binom{k-1}{i-1}\th^{k-i}{\om}(1-\th/2k)^t
+\frac{k^{2k}}{\th^i}\om(1-\th/3k)^t\binom{k}{2}\bfrac{\th}{k^2}^{i+1}(1+\th/k)\\
&\gap{2}+(1-\th(1-k^{2k}\e))k^{2k-2i}\om(1-\th/3k)^t+\D^{-\e/3}\\
&\leq
k^{2k-2i}(1-\th/3k)^{t+1}\om\brac{\binom{k-1}{i-1}\bfrac{\th}{k^2}^{k-i}\frac{(1-\th/2k)^t}{(1-\th/3k)^{t+1}}
+
\frac{\th(1+\th/k)+2(1-\th(1-k^{2k}\e))}{2(1-\th/3k)}}+\D^{-\e/3}\\
&\leq k^{2k-2i}(1-\th/3k
)^{t+1}\om\bfrac{2-\th(1-2k^{2k}\e)+O(\th^2)}{2(1-\th/3k)}+\D^{-\e/4}\\
&\leq k^{2k-2i}(1-\th/3k)^{t+1}\om.
\end{align*}
{\bf Property \eqref{C}:} Trivial.

{\bf Property \eqref{D}:}
If $d(u)\leq (1-\th/2k)^t\D$ then from
\eqref{D1} we get
\begin{eqnarray*}
d'(u)&\leq& \brac{1-\frac{\th}{k}}\brac{1-\frac{\th}{2k}}^{t}\D+\D^{2/3}\\
&=&\brac{1-\frac{\th}{2k}}^{t+1}\D-\frac{\th}{2k}\brac{1-\frac{\th}{2k}}^t\D+\D^{2/3}\\
&\leq& \brac{1-\frac{\th}{2k}}^{t+1}\D.
\end{eqnarray*}

{\bf Property \eqref{Dc}:}
$$d'_i(u,c)\leq (1+2k\th)^t\D\hp^{k-i}+2k\th(1+2k\th)^t\D\hp^{k-i}= (1+2k\th)^{t+1}\D\hp^{k-i}$$

To complete the proof it suffices to show that there are choices for $\g_u(c),\eta_u(c),\,u\in U,c\in C$
such that \eqref{A1}--\eqref{duc} hold.

In order to help understand the following computations, the reader
is reminded that quantities $\Xi_u,\F_{u,i},\omega,\th^{-1}$ can all
be upper bounded by $\D^{o(1)}$. Note also that in \eqref{Dc}, $(1+2k\th)^{t_0}=\log^{O(1)}\D=\D^{o(1)}$.
\subsection{Bad colors}\label{s1}
We now put a bound on the weight of the colors in $B(u)$.

Assume that \eqref{A}--\eqref{D} hold. It follows from \eqref{C}
that
\beq{n1}
h_u^{(0)}-h_u^{(t)}\leq k^{2k}\e
\sum_{i=0}^\infty(1-\th/3k)^i=3k\times k^{2k}{\om}=
\frac{3\e^2\log\D}{100}.
\eeq
Since $p_u^{(0)}(c)=1/q$ for all $u,c$
we have
\begin{eqnarray*}
h_u^{(0)}&=&-\sum_cp_u^{(0)}(c)\log p_u^{(0)}(c)\nonumber\\
&=&-\sum_cp_u^{(t)}(c)\log p_u^{(0)}(c)-(\log 1/q)\sum_c(p_u^{(0)}(c)-p_u^{(t)}(c))\\
&\geq &-\sum_cp_u^{(t)}(c)\log p_u^{(0)}(c)-t\D^{-\e}\log\D.
\end{eqnarray*}
where the last inequality uses \eqref{A}.

Plugging this lower bound on $\tz{h}_u$ into \eqref{n1} gives
\begin{eqnarray*}
\frac{3\e^2\log\D}{100}&\geq&h_u^{(0)}-\ts{h}_u\\
&\geq&-\sum_c\ts{p}_u(c)\log p_u^{(0)}(c)-t\D^{-\e}\log\D+\sum_c\ts{p}_u(c)\log\ts{p}_u(c)\\
&=&\sum_c\ts{p}_u(c)\log(\ts{p}_u(c)/p_u^{(0)}(c))-t\D^{-\e}\log\D.
\end{eqnarray*}
Thus,
\beq{S1}
\sum_cp_u^{(t)}(c)\log(p_u^{(t)}(c)/p_u^{(0)}(c))\leq
\frac{3\e^2\log\D}{100}+\D^{-\e/2}.
\eeq
Now, all terms in \eqref{S1} are
non-negative ($p_u^{(t)}(c)=0$ or $p_u^{(t)}(c)\geq
p_u^{(0)}(c)$). Thus after dropping the contributions from $c\notin B(u)$ we get
\begin{eqnarray*}
\frac{3\e^2\log\D}{100}+\D^{-\e/2}&\geq& \sum_{c\in B(u)}p_u^{(t)}(c)\log(p_u^{(t)}(c)/p_u^{(0)}(c))\\
&=&\sum_{c\in B(u)}p_u^{(t)}(c)\log(\hp q) =\sum_{c\in B(u)}p_u^{(t)}(c)\log(\D^{\e-o(1)})\\
&\geq& \frac23\e p_u(B(u))\log\D.
\end{eqnarray*}
 So,
\beq{star}
p_u(B(u))\leq \frac{\e}{10}.
\eeq

\subsection{Verification of Dynamics}\label{dyn}
Let
$\cE_{\ref{A1}}(u)$ -- $\cE_{\ref{duc}}(u)$ be the events claimed in
equations \eqref{A1} -- \eqref{duc}. Let
$\cE(u)=\cE_{\ref{A1}}(u)\cap \cdots\cap \cE_{\ref{duc}}(u)$. We
have to show that $\bigcap_{u\in U}\cE(u)$ has positive probability.
We use the local lemma. Each of the above events depends only on the
vertex $u$ or  its neighbors. Therefore, the dependency graph of the
$\cE(u),\,u\in U$ has maximum degree $\D^{O(1)}$ and so it is enough to
show that each event
$\cE_{\ref{A1}}(u),\ldots,\cE_{\ref{duc}}(u),\,u\in U$ has failure
probability $e^{-\D^{\O1}}$.

While parameters $\Xi_u,\F_u$ etc. are only needed for $u\in U$ we do not for example consider $\Xi_u'$ conditional on
$u\in U'$. We do not impose this conditioning and so we do not have to deal with it.
Thus the local lemma will guarantee a value for $\Xi_u,\,u\in U\setminus U'$
and we are free to disregard it for the next round. (We will however face this conditioning for other reasons, see \eqref{bugger}).

In the following we will use two forms of Hoeffding's inequality for
sums of bounded random variables: Suppose  first that
$X_1,X_2,\ldots,$ $X_m$ are independent random variables and
$|X_i|\leq a_i$ for $1\leq i\leq m$. Let $X=X_1+X_2+\cdots+X_m$.
Then, for any $t>0$,
\beq{hoef}
\max\set{\Pr(X-\E(X)\geq
t),\Pr(X-\E(X)\leq -t)}\leq \exp\set{-\frac{2t^2}{\sum_{i=1}^m
a_i^2}}.
\eeq
We will also need the following version in the special
case that $X_1,X_2,\ldots,X_m$ are independent [0,1] random variables.
For $\a>1$ we have
\beq{hoef1}
\Pr(X\geq L)\leq (3/\a)^{L}
\eeq
for
any $L\geq \a\E(X)$. (We replace $e$ by 3 as the symbol $e$ is over-used in the paper).

For proofs, see for example Alon and Spencer \cite{AS}, Appendix A and Lugosi \cite{Lugosi}.
\subsubsection{Dependencies}
In our random experiment, we start with
the $p_u(c)$'s and then we instantiate the independent random
variables $\g_u(c),\eta_u(c),u\in U,\,c\in C$ and then we compute
the $p_u'(c)$ from these values. Observe first that $p_u'(c)$
depends only on $\g_v(c),\eta_v(c)$ for $v=u$ or $v$ a neighbor of
$u$ in $H$. So $p_u'(c)$ and $p_v'(c^*)$ are independent if $c\neq
c^*$, even if $u=v$. We call this {\em color independence}.

Let
$$N_i(u)=\set{\set{u_2,u_3,\ldots,u_i}\subseteq U:\exists e\in H\ s.t.\ \set{u,u_2,\ldots,u_i}\subseteq e}.$$
We shorten $N_2(u)$ to $N(u)$.

If $f=\set{u_2,u_3,\ldots,u_i}\in N_i(u)$ and $l>i$, then
$E_{u,f,l}=\set{e\in N_l(u):e\supseteq f\cup\set{u},e\in H_l}$.

Next let
$$N_{i,l}(u)=\set{f=\set{u_2,\ldots,u_i}\subseteq U:\;f\in N_i(u)\ and\ E_{u,f,l}\neq\emptyset},\quad 2\leq i< l\leq k.$$
In words, $N_{i,l}$ is the collection of $i$-sets containing $u$
that are subsets of edges of $H_l$.

In these definitions $H_k=\ts{H}$.

For each $v\in N(u)$ we let
$$C_u(v)=\set{c\in C:\g_u(c)=1}\cup  L(v) \cup B(v).$$
Note that while the first two sets in this union depend on the random choices made in this
round, the set $B(v)$ is already defined at the beginning of the round.

We will later use the fact that if $c^*\notin C_u(v)$ and
$\g_v(c^*)=1$ then this is enough to place $c^*$ into $\Psi(v)$ and
allow $v$ to be colored. Indeed, we only need to check that $c^*
\not\in A^{(t-1)}(v)$ as it will then follow that $c^* \not\in
A(v)$. However, $\g_v(c^*)=1$ implies that $p_v(c^*) \ne 0$ from
which it follows that $c^* \not\in A^{(t-1)}(v)$.

Let $Y_v=\sum_cp_v(c)1_{c\in C_u(v)}=p_v(C_u(v))$. $C_u(v)$ is a
random set and $Y_v$ is the sum of $q$ independent random
variables each one bounded by $\hp$. Then by \eqref{N1}, \eqref{N2} and \eqref{crude1},
\begin{eqnarray*}
\E(Y_v)&\leq&\sum_{c\in C}p_v(c)\Pr(\g_u(c)=1)+\sum_{c\in C}p_v(c)(1-q_v(c))+p_v(B(v))\\
&\leq&\th\sum_{c\in C}p_u(c)p_v(c)+\th^{k-1}\Xi_v+\sum_{i=2}^{k-1}\th^{i-1} \F_{v,i} +p_v(B(v)).
\end{eqnarray*}
Now let us bound each term separately:
$$\th\sum_{c\in C}p_u(c)p_v(c)\le \th q \hp^2 <\th
\D^{1/(k-1)}\D^{2\e-2/(k-1)}<\frac{\e}{3}.$$
Using (\ref{B}) we obtain
$$\th^{k-1}\Xi_v<\om\th^{k-1}+t\th^{k-1}\D^{-\e}\le
\e\th^{k-2}+t\th^{k-1}\D^{-\e}<\frac{\e}{6}+\frac{\e}{6}=\frac{\e}{3}.$$
Using (\ref{Ca}) we obtain
$$\sum_{i=2}^{k-1}\th^{i-1} \F_{v,i}\le \sum_{i=2}^{k-1}\th^{i-1}k^{2k-2i}\om(1-\th/3k)^t\leq k^{2k-1}\e.$$
Together with $\Pr(B(v)) \le \e/10$ we get
$$\E(Y_v) \le (k^{2k-1}+1)\e.$$
Hoeffding's inequality then gives
$$\Pr(Y_v\geq \E(Y_v)+\r)\leq
\exp\set{-\frac{2\r^2}{q\hp^2}}<e^{-2\r^2\D^{1/k}}.$$
Taking
$\r=\D^{-1/2k}$ say, it follows that
\beq{Cv1}
\Pr(p_v(C_u(v))\geq
(k^{2k-1}+2)\e)=\Pr(Y_v\geq (k^{2k-1}+2)\e)\le e^{-\D^{1/2k}}.
\eeq
Let $\cE_{\eqref{Cv1}}$ be the event $\set{p_v(C_u(v))\leq
(k^{2k-1}+2)\e}$.

Now consider some fixed vertex $u\in U$. It will sometimes be
convenient to condition on the values $\g_x(c),\eta_x(c)$ for all
$c\in C$ and all $x\notin N(u)$
 and for $x=u$.
This conditioning is needed to obtain independence.
We let $\cC$ denote these conditional values.

\begin{remark} \label{rem} Note that $\cC$ determines  the set $C_u(v)$, and hence it also
determines whether or not $\cE_{\eqref{Cv1}}$ occurs. Indeed, if
$\g_u(c)=1$, then $c \in C_u(v)$.  On the other hand,  if
$\g_u(c)=0$ then whether or not $c\in L(v)$ depends only on colors
tentatively assigned to vertices not in $N(u)$. This uses the
simplicity and triangle-freeness of $H$.
\end{remark}

Given the conditioning $\cC$, simplicity and triangle freeness imply
that the events $\set{v\notin U'}$, $\set{w\notin U'}$ for $v,w\in
N(u)$ are independent provided $u,v,w$ are not part of an edge of
$H$. Indeed, triangle-freeness implies that in this case, there is
no edge containing both $v$ and $w$.  Therefore the random choices
at $w$ will not affect the coloring of $v$ (and vice versa). Thus
random variables $p'_{v}(c),p'_{w}(c)$ will become (conditionally)
independent under these circumstances. We call this {\em conditional
neighborhood independence}.

{\bf 3.9.1.1\ \ Some expectations}

Let us fix a color $c$ and an edge $u_1u_2\cdots u_k\in H$ (here we mean $H$
and not $H^{(t)}$). In this subsection we will
estimate the expectation of $\E\brac{\prod_{j=1}^ip'_{u_j}(c)}$ for $2\leq i\le k$ in two distinct situations.

{\bf Case 1:} $i=k$ and $u_1u_2\cdots u_k\in \ts{H}$.

Our goal is to prove that
\beq{eprod}
\E\brac{\prod_{j=1}^kp'_{u_j}(c)}\leq
(1+2k\th^2\hp^2)\prod_{j=1}^kp_{u_j}(c)
\eeq
If $c\in
\bigcup_{j=1}^kA^{(t-1)}(u_j)$ then $\prod_{j=1}^kp'_{u_j}(c)=0$.
Assume then that $c\notin \bigcup_{j=1}^kA^{(t-1)}(u_j)$. If  Case B
of \eqref{pdash} occurs for any of $u_1,u_2,\ldots,u_k$ e.g. $u_k$
then
\beq{factorout}
\E\brac{\prod_{j=1}^kp'_{u_j}(c)\biggr| Case\ B\
for\ k}=\E\brac{\prod_{j=1}^{k-1}p'_{u_j}(c)}p_{u_k}(c).
\eeq
This
is because in Case B the value of $\eta_{u_k}(c)$ is independent of
all other random variables and we may use \eqref{puc1}. One can see
then that we have to prove something slightly more general than
\eqref{eprod}. So we now aim to show that
\beq{eprod1}
\E\brac{\prod_{j=1}^ip'_{u_j}(c)}\leq
(1+2k\th^2\hp^2)\prod_{j=1}^ip_{u_j}(c)
\eeq
assuming that $1\leq
i\leq k$ and that there is an edge $u_1u_2\cdots u_k\in \ts{H}$ and
that Case A of \eqref{pdash} happens for $u_j,c,\,1\leq j\leq i$.
The case $i=1$ follows from \eqref{puc1} and so we assume that
$i\geq 2$.  By simplicity, $u_{i+1},\ldots,u_k$ are determined by
$u_1,u_2,\ldots,u_i$.

Now $\prod_{j=1}^ip'_{u_j}(c)=0$ unless $c\notin \bigcup_{j=1}^iL(u_j)$. Consequently,
\beq{euvd}
\E\brac{\prod_{j=1}^ip'_{u_j}(c)}=\prod_{j=1}^i\frac{p_{u_j}(c)}{q_{u_j}(c)}\times
\Pr\brac{c\notin \bigcup_{j=1}^iL(u_j)}.
\eeq
Furthermore,
\begin{eqnarray*}
\Pr\brac{c\notin \bigcup_{j=1}^iL(u_j)\biggr|
\g_{u_j}(c)=0,1\leq j\leq k}&=&\prod_{j=1}^i\brac{q_{u_j}(c)\brac{1-\th^{k-1}\prod_{j'\neq j}p_{u_{j'}}(c)}^{-1}}\\
&\leq&(1+2k\th^{k-1}\hp^{k-1})\prod_{j=1}^iq_{u_j}(c)\\
&\leq&(1+2k\th^2\hp^2)\prod_{j=1}^iq_{u_j}(c).
\end{eqnarray*}
On the other hand we will show that
\beq{gamm01}
\Pr\brac{c\notin
\bigcup_{j=1}^iL(u_j)\biggr|\exists 1\leq j\leq k:\g_{u_j}(c)=1}\leq
\Pr\brac{c\notin \bigcup_{j=1}^iL(u_j)\biggr| \g_{u_j}(c)=0,1\leq
j\leq k}.
\eeq
This is intuitively clear, since color $c$ is at
least as likely to be lost at $u_j$ if it is tentatively activated
at some $u_{j'}$. Indeed, to see this formally, partition the
probability space $\Omega$ of outcomes of the $\g$'s and $\eta$'s
into the sets $\Omega_{\e_{1},\ldots,\e_k}$ in which
$\g_{u_j}(c)=\e_j\in \set{0,1}$ for $1\leq j\leq k$. Let
$\Omega_{\e_{1},\ldots,\e_k}'$ be the set of outcomes in
$\Omega_{\e_{1},\ldots,\e_k}$ in which $c\notin
\bigcup_{j=1}^iL(u_j)$. Now consider the map
$f:\Omega_{\e_{1},\ldots,\e_k}\to \Omega_{0,\ldots,0}'$ which just
sets $\g_{u_j}(c)$ to 0 for $1\leq j\leq k$. Then if $\p_j^1=\th
p_{u_j}(c)$ and $\p_j^0=1-\p_j^1$
$$\frac{\Pr(\Omega_{\e_{1},\ldots,\e_k})}{\Pr(\Omega_{0,\ldots,0})}=
\frac{\prod_{j=1}^k \p_j^{\e_j}}{\prod_{j=1}^k\p_j^0}=\frac{\Pr(\Omega'_{\e_{1},\ldots,\e_k})}
{\Pr(f(\Omega'_{\e_{1},\ldots,\e_k}))}.$$
If $c \not \in \bigcup_{j=1}^iL(u_j)$ and $\exists 1\leq j\leq k$ such that $\gamma_{u_j}(c)=1$, then we still
have $c \not\in \bigcup_{j=1}^iL(u_j)$ if we change $\gamma_{u_j}(c)$ to 0 for $1\leq j\leq k$ and
make no other changes. Consequently, $f(\Omega_{\e_{1},\ldots,\e_k}')\subseteq
\Omega_{0,\ldots,0}'$ and we have
$$\frac{\Pr(\Omega_{0,\ldots,0}')}{\Pr(\Omega_{0,\ldots,0})}
\geq \frac{\Pr(f(\Omega_{\e_{1},\ldots,\e_k}'))}{\Pr(\Omega_{\e_{1},\ldots,\e_k})}\cdot
\frac{\Pr(\Omega_{\e_{1},\ldots,\e_k})}{\Pr(\Omega_{0,\ldots,0})}
=\frac{\Pr(\Omega_{\e_{1},\ldots,\e_k}')}{\Pr(\Omega_{\e_{1},\ldots,\e_k})},$$
which is \eqref{gamm01}.

It follows that
\begin{equation} \label{starstar}
\Pr\brac{c\notin \bigcup_{j=1}^iL(u_j)}\leq
(1+2k\th^2\hp^2)\prod_{j=1}^iq_{u_j}(c).
\end{equation}
 and in combination with \eqref{euvd} this proves \eqref{eprod1} and hence \eqref{eprod}.

{\bf Case 2:} $e=u_1u_2\cdots u_i\in H_i,\,\k(e)=c$.

Our goal is now to prove
\beq{bugger}
\E\brac{\prod_{j=1}^ip_{u_j}'(c)\times 1_{u_1,u_2,\ldots,u_i\in U'}}\leq (1+k\th\hp)
\prod_{j=1}^ip_{u_j}(c).
\eeq

Suppose that $p'_{u_j}(c)$ is determined by Case A of \eqref{pdash} for $1\leq j\leq l$ and by Case B otherwise.
We factor out $\prod_{j=l+1}^ip_{u_j}(c)$ as in \eqref{factorout} and concentrate on  bounding
\begin{align}
&\E\brac{ \prod_{j=1}^lp_{u_j}'(c)\times 1_{u_1,u_2,\ldots,u_i\in U'}}\nonumber\\
&=\prod_{j=1}^l\frac{p_{u_j}(c)}{q_{u_j}(c)} \times\Pr(c\notin L(u_1)\cup \cdots L(u_l)
\wedge u_1,\ldots,u_i\in U')\nonumber\\
&\leq\prod_{j=1}^l\frac{p_{u_j}(c)}{q_{u_j}(c)} \times\Pr(c\notin L(u_1)\cup \cdots L(u_l))\nonumber\\
&\leq\prod_{j=1}^l\frac{p_{u_j}(c)}{q_{u_j}(c)} \times\Pr(c\notin L(u_1)\cup \cdots L(u_l)\mid \g_{u_1}(c)=\cdots=
\g_{u_l}(c)=0)\label{nonumber}\\
&\leq \prod_{j=1}^l p_{u_j}(c)\times (1-\th\hp)^{-l}.\label{pol}\\
&\leq (1+k\th\hp) \prod_{j=1}^l p_{u_j}(c)
\end{align}
and \eqref{bugger} follows.

{\bf Explanation:} Equation \eqref{nonumber} follows as for \eqref{gamm01}. Equation \eqref{pol} now follows because
the events $c\notin L(u_j)$ become conditionally independent. And then $\Pr(c\notin L(u_j)\mid \g_{u_j}(c)=0)$
gains a factor $\brac{1-\th^{i-1}\prod_{j'\neq j} p_{u_j}(c)}^{-1}\leq (1-\th\hp)^{-1}$.
\subsubsection{Proof of \eqref{A1}}
Given the $p_u(c)$ we see that
if $Z'=\sum_{c\in C}p_u'(c)$ then $\E(Z')=\sum_{c\in C}p_u(c)$. This
follows on using \eqref{puc1}. By color independence $Z'$ is the sum
of $q$ independent non-negative random variables each bounded by
$\hp$. Applying \eqref{hoef} we see that
$$\Pr(|Z'-\E(Z')|\geq \r)\leq 2\exp\set{-\frac{2\r^2}{q\hp^2}}=2e^{-2\r^2\D^{1/(k-1)-2\e-o(1)}}.$$
We take $\r=\D^{-\e}$ to see that $\cE_{\ref{A1}}(u)$ holds
\whp\footnote{By \whp, {\em with high probability}, we mean with
probability $1-e^{-\D^{\Omega(1)}}$.}.
\subsubsection{Proof of \eqref{B1}}
Given the $p_u(c)$ we see that by \eqref{eprod},
$\Xi_{e}'$ has expectation no more than $\Xi_{e}(1+2k\th^2\hp^2)$ and
is the sum of $q$ independent non-negative random variables, each of
which is bounded by $\hp^k$. We have used color independence again
here. Applying \eqref{hoef} we see that
$$\Pr(\Xi_{e}'\geq \Xi_e(1+2k\th^2\hp^k)+\r/2)\leq \exp\set{-\frac{\r^2}{2q\hp^{2k}}}\leq
 e^{-\r^2\D^{(2k-1)/(k-1)-2k\e-o(1)}}.$$
We also have
$$k\Xi_{e}\th^2 \hp^2\le
k\left(\frac{\om}{\D}+
\frac{t}{\D^{1+\e}}\right)\th^2\hp^2<\frac{1}{2\D^{1+\e}}.$$
We take $\r=\D^{-1-\e}$ to obtain
$$\Pr(\Xi_e'\geq \Xi_e+\D^{-1-\e})\leq e^{-\D^{\Omega(1)}}$$
and so $\cE_{\ref{B1}}(u)$ holds \whp.

\subsubsection{Proof of \eqref{B1a}}

Throughout this section $u_1=u$. Recall that
$$\F_{u,i}=\sum_{c\in C}\ \sum_{e=uu_2\cdots u_i\in H_i}1_{\k(e)=c}\prod_{j=1}^ip_{u_j}(c).$$
If $\set{u_2,\ldots, ,u_i}\in N_i(u)$ and $e=uu_2\cdots u_i\notin
H_i$ then $\k(e)$ is defined to be $0\notin C$. Now
$$
\F_{u,i}'-\F_{u,i}=\sum_{c\in C}\brac{\sum_{e=uu_2\cdots u_i \in
H_i'}1_{\k'(e)=c}\,\prod_{j=1}^ip_{u_j}'(c) -\sum_{e=uu_2\cdots
u_i\in H_i}1_{\k(e)=c}\,\prod_{j=1}^ip_{u_j}(c)}$$
 If $e\in H_i$ and
$\k(e)=c$ then $e\in H_i'$ is equivalent to $\k'(e)=c$. If $\k(e)\ne
c$ but $\k'(e)=c$ then the edge of $H$ containing $u,u_2, \ldots,
u_i$ has some other vertices in $U$ that will be colored with $c$ in
the current round. Consequently, the above expression is
$$D_1+\sum_{l=i+1}^kD_{2,l}$$
where
\begin{eqnarray*}
D_1&=&\sum_{c\in C}\sum_{\substack{e=uu_2\cdots u_i\in
 H_i\\ \k(e)=c}} \brac{1_{\k'(e)=c}\,\prod_{j=1}^ip_{u_j}'(c)-\prod_{j=1}^ip_{u_j}(c)}\\
D_{2,k}&=&\sum_{c\in C}\sum_{\set{u_2,\ldots u_i}\in N_{i,k}(u)}
 1_{\k'(uu_2\cdots u_i)=c}\,\prod_{j=1}^ip_{u_j}'(c)\\
D_{2,l}&=&\sum_{c\in C}\sum_{\substack{\set{u_2,\ldots u_i}\in
N_{i,l}(u)\\ \k(uu_2\cdots u_l)=c}}
 1_{\k'(uu_2\cdots u_i)=c}\,\prod_{j=1}^ip_{u_j}'(c)\qquad i+1\leq l\leq k-1.
\end{eqnarray*}

Here $D_1$ accounts for the contribution from edges leaving $H_i$
and $D_{2,i+1},\ldots,D_{2,k}$ account for the contribution from
edges entering $H_i$.

We bound $\E(D_1)$ and $\E(D_{2,i+1}),\ldots,\E(D_{2,k})$ separately.

$\E(D_1)$:\\
\begin{eqnarray*}
D_1&=&\sum_{c\in C}\sum_{\substack{e=uu_2\cdots u_i\in
 H_i\\ \k(e)=c}} \brac{1_{\k'(e)=c}\,\prod_{j=1}^ip_{u_j}'(c)-\prod_{j=1}^ip_{u_j}(c)}\\
&=&-\sum_{c\in C}\sum_{\substack{e=uu_2\cdots u_i\in H_i\\ \k(e)=c \\ \k'(e)\ne c}}\prod_{j=1}^ip_{u_j}(c)
+\sum_{c\in C}\sum_{\substack{e=uu_2\cdots u_i\in H_i\\ \k(e)=c \\ \k'(e)=c}}
\brac{\prod_{j=1}^ip_{u_j}'(c)- \prod_{j=1}^ip_{u_j}(c)}.\\
\end{eqnarray*}
Now suppose that $\exists j\ge2: u_j \not\in U'$. This means that $u_j$ has been colored in the current round and
so $e \not\in H_i'$.  So $\k'(e) \ne c$ is implied by $\exists j\ge2: u_j \not\in U'$ and more simply
it is implied by $u_2\notin U'$. Conversely, if $e\in H_i'$ then $u_1,u_2,\ldots,u_i\in U'$.
Therefore the prior expression is bounded from above by
$$-D_{1,1}+D_{1,2}$$
\begin{eqnarray*}
\noalign{where}\\
D_{1,1}&=&\sum_{c\in C}\sum_{\substack{e=uu_2\cdots u_i\in H_i\\ \k(e)=c}}\prod_{j=1}^ip_{u_j}(c)1_{u_2
\notin U'}\\
D_{1,2}&=&\sum_{c\in C}\sum_{\substack{e=uu_2\cdots u_i\in H_i\\ \k(e)=c}}
\brac{\brac{\prod_{j=1}^ip_{u_j}'(c)-\prod_{j=1}^ip_{u_j}(c)}\times 1_{u_1,u_2,\ldots,u_i\in U'}}.\\
\end{eqnarray*}

Suppose that $e=uu_2\cdots u_k\in H$ and $u_{i+1},u_{i+2},\ldots,u_k\notin U$ and $\k(u_{i+1})=\cdots=\k(u_k)=c$.
Let $v=u_2$. Recall that
$$C_u(v)=\set{c\in C:\g_u(c)=1}\cup  L(v) \cup B(v).$$
If there is a tentatively activated color $c^*$ at $v$ (i.e.
$\gamma_v(c^*)=1$) that lies outside $C_u(v) \cup \{c\}$, then $c^*
\in \Psi(v)$ and $v$ will be colored in this round (recall that we
had argued earlier that $c^* \not\in A^{(t-1)}(v)$). Therefore
$$\Pr(v\notin U'\mid\cC)\geq \Pr(\exists c^*\notin C_u(v)\cup\set{c}:\; \g_v(c^*)=1\mid\cC).$$
We have introduced the conditioning $\cC$ because we will need it later when we prove concentration.

So by inclusion-exclusion and the independence of the $\g_v(c^*)$ we can write
\begin{eqnarray*}
\E\brac{1_{v\notin U'}\mid\cC}&\geq&\Pr(\exists c^*\notin C_u(v)\cup\set{c}:\; \g_v(c^*)=1\mid\cC)\\
&\geq&\sum_{c^*\notin C_u(v)\cup\set{c}}\Pr(\g_v(c^*)=1\mid\cC)-\frac12\sum_{c_1^*\neq c_2^*\notin C_u(v)\cup\set{c}}\Pr(\g_v(c_1^*)=\g_v(c_2^*)=1\mid\cC)\\
&\geq& \sum_{c^*\notin C_u(v)\cup\set{c}}\th
p_v(c^*)-\frac12\brac{\sum_{c^*\notin C_u(v)\cup\set{c}}
\th p_v(c^*)}^2 \end{eqnarray*}
Now
\begin{eqnarray*}
\sum_{c^*\notin C_u(v)\cup\set{c}}\th
p_v(c^*)
&=&\sum_{c^* \in C}\th
p_v(c^*)- \sum_{c^* \in C_u(v)} \th p_v(c^*) -\th p_v(c)\\
&\ge& \th((1-t\D^{-\e}) -p_v(C_u(v)) -\hp) \\
&>& \th(1-p_v(C_u(v)) -\e/2)
\end{eqnarray*}
where we have used (\ref{A}).
Also by (\ref{A}) and the definition of $\hp$ we have
$$\sum_{c \ne c^*} p_v(c^*) \le 1+(t+1)\D^{-\e} <1.1.$$
Consequently
$$\frac12\brac{\sum_{c^*\notin C_u(v)\cup\set{c}}
\th p_v(c^*)}^2=\frac{\th^2}{2}\brac{\sum_{c^*\notin C_u(v)\cup\set{c}}
p_v(c^*)}^2 \le \frac{2\th^2}{3} <\frac{\th \e}{2}.$$
Putting these facts together yields
$$\E\brac{1_{v\notin U'}\mid\cC} \ge \th(1-p_v(C_u(v))-\e).$$
Therefore
$$
\E(D_{1,1}\mid\cC)\geq \th(1-p_v(C_u(v))-\e)
\sum_{c\in C}\sum_{\substack{e=uu_2\cdots u_i\in H_i\\ \k(e)=c}}\prod_{j=1}^ip_{u_j}(c)=
\th(1-p_v(C_u(v))-\e)\F_{u,i}.
$$
Given $\cC$, Remark \ref{rem} implies that  $\cE_{\eqref{Cv1}}$
either holds for all outcomes in $\cC$, or fails for all outcomes in
$\cC$.  So,
 \beq{D11}
\E(D_{1,1}\mid\cC)\geq \th(1-(k^{2k-1}+3)\e)\F_{u,i},\qquad for\
\cC\ such\ that\ \cE_{\eqref{Cv1}}\ occurs. \eeq We now consider
$D_{1,2}$. It follows from \eqref{Ca} that $\F_{u,i} <k^{2k-2i}\om$.
Together with \eqref{bugger}, this gives \beq{D12} \E(D_{1,2})\leq
k\F_{u,i}\th\hp\leq k^{2k}\e\hp. \eeq

$\E(D_{2,k})$:\\
Recall that
$$D_{2,k}=\sum_{c\in C}\sum_{\{u_2, \ldots, u_i\}\in N_{i,k}(u)}1_{\k'(uu_2\cdots u_i)=c}\prod_{j=1}^ip_{u_j}'(c).$$
Instead of summing over sets in $N_{i,k}(u)$, we may sum over edges
$uu_2\cdots u_k \in H^{(t)}$, and then over subsets of these edges
that lie in $N_{i,k}(u)$. Thus
$$D_{2,k}=\sum_{c\in C}\sum_{uu_2\cdots u_k\in \ts{H}}\sum_{f \subset
\{u_2, \ldots, u_k\}, |f|=i-1} 1_{\k'(f \cup \{u\})=c}\prod_{u_j\in
f \cup \{u_1\}}p_{u_j}'(c).$$

Fix an edge $uu_2\cdots u_k\in H^{(t)}$. If $u_{i+1},\ldots,u_k$ are
colored with $c$ in this round, then certainly $c$ must have been
tentatively activated at these vertices. Therefore,
\begin{eqnarray}
\E\brac{1_{\k'(uu_2\cdots u_i)=c}\prod_{j=1}^ip_{u_j}'(c)}&\leq&\E\brac{\prod_{j=i+1}^k\g_{u_j}(c)
\prod_{j=1}^ip_{u_j}'(c)}\nonumber\\
&\leq&\th^{k-i} \prod_{j=i+1}^kp_{u_j}(c)\prod_{j=1}^i\frac{p_{u_j}(c)}{q_{u_j}(c)}
\Pr\brac{c\notin \bigcup_{j=1}^iL(u_j)\biggr| \bigwedge_{j=i+1}^k(\g_{u_j}(c)=1)}\nonumber\\
&\leq&\th^{k-i} \prod_{j=i+1}^kp_{u_j}(c)\prod_{j=1}^i\frac{p_{u_j}(c)}{q_{u_j}(c)}
\Pr\brac{c\notin \bigcup_{j=1}^iL(u_j)}\label{imitrex}\\
&\leq&\th^{k-i}\prod_{j=1}^kp_{u_j}(c) (1+2k\th^2\hp^2).\label{imit}
\end{eqnarray}
We use the argument for \eqref{gamm01} to obtain \eqref{imitrex} and (\ref{starstar}) to obtain (\ref{imit}).

It follows that
\beq{ED2k}
\E(D_{2,k})\leq \binom{k-1}{i-1}\th^{k-i} \Xi_u(1+2k\th^2\hp^2).
\eeq
$\E(D_{2,l}),\,l<k$:\\
Recall that
$$D_{2,l}=\sum_{c\in C}\sum_{\substack{\set{u_2,u_3,\ldots u_i}\in N_{i,l}(u)\\ \k(uu_2\cdots u_l)=c}}
 1_{\k'(uu_2\cdots u_i)=c}\,\prod_{j=1}^ip_{u_j}'(c).$$
Fix an edge $uu_2\cdots u_k\in H$ with $uu_2\cdots u_l\in H_l$. Then arguing as we did for \eqref{imit} we have
$$\E\brac{1_{\k'(uu_2\cdots u_i)=c}\prod_{j=1}^ip_{u_j}'(c)}\leq\th^{l-i}\prod_{j=1}^lp_{u_j}(c) (1+2k\th^2\hp^2).$$
It follows that \beq{ED2l} \E(D_{2,l})\leq \binom{l-1}{i-1}\th^{l-i}
\F_{u,l}(1+2k\th^2\hp^2). \eeq {\bf 3.9.4.1\ \ Concentration}

We first deal with $D_{1,1}$. For this we condition on the values
$\g_w(c),\eta_w(c)$ for all $c\in C$ and all $w\notin N(u)$ and for
$w=u$. Then by conditional neighborhood independence $D_{1,1}$ is
the sum of at most $d_i(u)$ independent random variables of value at
most $\hp^i$. By (\ref{Dc}), we have $d_i(u) \le q(1+2k\th)^t\D\hp^{k-i}=\D^{1+1/(k-1)+o(1)}\hp^{k-i}$. So,
for $\r>0$,
$$\Pr(D_{1,1}-\E(D_{1,1}\mid \cC)\leq -\r\mid\cC)\leq
\exp\set{-\frac{2\r^2}{d_i(u)\hp^{2i}}}\le
\exp\set{-\frac{2\r^2}{\D^{1+1/(k-1)+o(1)}\hp^{k-i}\hp^{2i}}}\le
e^{-\r^2\D^{i/k}}.$$

So, by \eqref{D11},
\begin{align}
&\Pr(D_{1,1}\leq \th(1-(k^{2k-1}+3)\e)\F_{u,i}-\D^{-1/2k})\nonumber\\
&=\sum_{\cC}\Pr(D_{1,1}\leq
\th(1-(k^{2k-1}+3)\e)\F_{u,i}-\D^{-1/2k}\mid\cC)\Pr(\cC)\nonumber\\
&\leq\sum_{\cC:\cE_{\eqref{Cv1}}\ occurs}\Pr(D_{1,1}\leq
\th(1-(k^{2k-1}+3)\e)\F_{u,i}-\D^{-1/2k}\mid\cC)\Pr(\cC)+
\Pr(\neg\cE_{\eqref{Cv1}})\nonumber\\
&\leq\sum_{\cC:\cE_{\eqref{Cv1}}\ occurs}\Pr(D_{1,1}\leq
\E(D_{1,1}\mid \cC)-\D^{-1/2k}\mid\cC)\Pr(\cC)+
\Pr(\neg\cE_{\eqref{Cv1}})\nonumber\\
&\leq e^{-\D^{1/k}}+e^{-\D^{1/2k}}\nonumber\\
&=e^{-\D^{\Omega(1)}}.\label{D11a}
\end{align}
Now consider the sum $D_{1,2}$. Let
$$Y_c=\sum_{\substack{e=uu_2\cdots u_i\in H_i\\ \k(e)=c\\\k'(e)=c}}
\brac{\prod_{j=1}^ip_{u_j}'(c)-\prod_{j=1}^ip_{u_j}(c)}.$$
$D_{1,2}$
is the sum of $q$ independent random variables $Y_c$ satisfying
$0\leq Y_c\leq d_c\hp^i$ where $d_c=d_i(u,c)$. Note that \eqref{Dc}
implies $d_c\leq \D^{1+o(1)}\hp^{k-i}$.

So, for $\r>0$,
$$\Pr(D_{1,2}-\E(D_{1,2})\geq \r)\leq \exp\set{-\frac{2\r^2}{\sum_cd_c^2\hp^{2i}}}\leq
\exp\set{-\frac{2\r^2}{\D^{2+1/(k-1)+o(1)}\hp^{2k}}} \leq
e^{-\r^2\D^{1/(k-1)-2k\e+o(1)}}.$$
We take $\r=\D^{-1/2k}$ to see
that $\Pr(D_{1,2}\geq 2\D^{-1/2k})\leq e^{-\D^\e}$. Combining this
with \eqref{D11a} we see that
\begin{align}
&\Pr(D_1\geq -\th(1-(k^{2k-1}+3)\e)\F_{u,i}+3\D^{-1/2k})\nonumber\\
&\leq \Pr(D_{1,1}\leq \th(1-(k^{2k-1}+3)\e)\F_{u,i}+\D^{-1/2k}) +
\Pr(D_{1,2}\geq
2\D^{-1/2k})\nonumber\\
&\leq e^{-\D^{\Omega(1)}}.\label{D11b}
\end{align}

We now deal with the $D_{2,l}$. There is a minor problem in that the
$D_{2,l}$ are sums of random variables for which we do not have a
sufficiently small absolute bound. These variables do however have a
small bound which holds with high probability. There are several
ways to use this fact. We proceed as follows: First assume $l\leq
k-1$ and let
$$D_{2,l,c}=\sum_{\substack{\set{u_2,u_3,\ldots u_i}\in N_{i,l}(u)\\ \k(uu_2\cdots u_l)=c}}
 1_{\k'(uu_2\cdots u_i)=c}\,\prod_{j=1}^ip_{u_j}'(c)$$
which we re-write as
$$D_{2,l,c}=\sum_{\substack{e=uu_2\cdots u_l\in H_l\\\k(e)=c}}Z_e,$$
where
$$Z_{uu_2\cdots u_l}=\sum_{\substack{S \subset\set{u_2,u_3,\ldots u_l}\\|S|=i-1}}1_{\k'(S \cup
\{u\})=c}\, \prod_{u_j\in S \cup \{u_1\}}p_{u_j}'(c).$$
Then we
let
$$\hD_{2,l}=\sum_{c\in C}\min\set{(1+2k\th)^t\D\hp^k,D_{2,l,c}}.$$
Observe that $\hD_{2,l}$ is the sum of $q$ independent random
variables each bounded by $(1+2k\th)^t\D\hp^k$. So, for $\r>0$,
$$\Pr(\hD_{2,l}-\E(\hD_{2,l})\geq \r)\leq \exp\set{-\frac{2\r^2}{\D^{2+o(1)}\hp^{2k}}}
\leq e^{-\r^2\D^{1/(k-1)-2k\e}}.$$
We take $\r=\D^{-1/2k}$ to see
that
\beq{hD2}
\Pr(\hD_{2,l}\geq \E(\hD_{2,l})+\D^{-1/2k})\leq
e^{-\D^{\e}}.
\eeq
We must of course compare $D_{2,l}$ and
$\hD_{2,l}$. Now $D_{2,l}\neq \hD_{2,l}$ only if there exists $c$
such that $D_{2,l,c}>(1+2k\th)^t\D\hp^k$. For each $c$, $D_{2,l,c}$ is the sum
of the $d_l(u,c)\leq (1+2k\th)^t\D\hp^{k-l}$ variables $Z_e,\,e\in H_l$. Each
$Z_e$ is bounded above by $\binom{l-1}{i-1} \hp^i$
and $\E(Z_e)\leq \binom{l-1}{i-1}\th^{l-i}\hp^l$. This is because $Z_e$
is bounded by the sum of $\binom{l-1}{i-1}$ variables $Z_{e,S}$, each
taking the value 0 or $\hp^i$. Here $Z_{e,S}$ corresponds to some
$S=\set{u_2,u_3,\ldots u_i} \subseteq \set{u_2,u_3,\ldots u_l}$.
Furthermore, $\Pr(Z_{e,S}=\hp^i)\leq (\th\hp)^{l-i}$ because this
will happen only if the vertices in $S$ tentatively choose $c$.

$H$ being simple and triangle free, if we condition on $\cC$ then the random variables $Z_e$ become independent.

Now put $X_e=Z_e/(\binom{l-1}{i-1}\hp^i)$ and $X=\sum_eX_e$. We see that
$0\leq X_e\leq 1$ and $\E(X_e)\leq (\th\hp)^{l-i}$.

We now use \eqref{hoef1} with $\alpha=1/(\binom{l-1}{i-1}\th^{l-i})$ and
$\E(X)\leq (1+2k\th)^t\D\hp^{k-l}\times (\th\hp)^{l-i}=
\frac{(1+2k\th)^t\D\hp^{k}}{\a\binom{l-1}{i-1}\hp^{i}}$. This gives
$$\Pr(D_{2,l,c}\geq (1+2k\th)^t\D\hp^k\mid\cC)\leq \Pr\brac{X\geq \frac{(1+2k\th)^t\D\hp^k}{\binom{l}{i}\hp^{i}}\biggr|\cC}
\leq \brac{3\binom{l}{i}\th^{l-i}}
^{(1+2k\th)^t\D\hp^{k-i}/\binom{l}{i}}$$
Therefore
\begin{eqnarray}
\Pr(D_{2,l}\neq \hD_{2,l})&\leq &\sum_\cC\Pr(\exists c:\;D_{2,l,c}\geq (1+2k\th)^t\D\hp^k\mid\cC)\Pr(\cC)\nonumber\\
&\leq&\sum_\cC q\brac{3\binom{l}{i}\th^{l-i}}^{(1+2k\th)^t\D\hp^{k-i}/\binom{l}{i}}\Pr(\cC)\nonumber\\
&\leq&e^{-\D^\e}.\label{hD2a}
\end{eqnarray}
It follows from \eqref{hD2a} and $\hD_{2,l}\leq D_{2,l}\leq \D$ that
$$|\E(D_{2,l})-\E(\hD_{2,l})|\leq \D\Pr(D_{2,l}\neq \hD_{2,l})
\leq\D qe^{-\D^\e}<\D^{-1/2k}.$$
Applying \eqref{hD2} and
\eqref{hD2a} we see that
\beq{D2l}
\Pr(D_{2,l}\geq \E(D_{2,l})+2\D^{-1/2k})\leq\\
\Pr(\hD_{2,l}\geq \E(\hD_{2,l})+\D^{-1/2k})+\Pr(D_{2,l}\neq
\hD_{2,l})\leq 2e^{-\D^{\e}}.
\eeq

We must now deal with the case of $l=k$ i.e.
$$D_{2,k,c}=\sum_{\set{u_2,u_3,\ldots u_i}\in N_{i,k}(u)}
 1_{\k'(uu_2\cdots u_i)=c}\,\prod_{j=1}^ip_{u_j}'(c)$$
and
$$\hD_{2,k}=\sum_{c\in C}\min\set{\D\hp^k,D_{2,k,c}}.$$
We re-write
$$D_{2,k,c}=\sum_{S\in N_{i,k}(u)}W_S$$
where for $S=\set{u_2,u_3,\ldots u_i}$,
$$W_S=\sum_{e\supseteq S,e\in \ts{H}}1_{\k'(uu_2\cdots u_i)=c}\prod_{j=1}^ip_{u_j}'(c).$$
Now we view $D_{2,k,c}$ as the sum of at most $\D$ random variables,
each of which is bounded by $\binom{k}{i}\hp^i$ and has expectation
bounded by $\binom{k}{i}\th^{k-i}\hp^k$. We now simply follow the
argument for $l<k$ by taking $l=k$ to show that
\beq{D2k}
\Pr(D_{2,k}\geq \E(D_{2,k})+2\D^{-1/2k})\leq 2e^{-\D^{\e}}.
\eeq
Indeed, \eqref{hD2} holds with $l=k$. Then
$$\Pr(D_{2,k}\neq \hD_{2,k})\leq q\Pr\brac{Bin\brac{\D,(\th\hp)^{k-i}}\geq \D\hp^{k-i}}
\leq q\brac{3\th^{k-i}}^{\D\hp^{k-i}}\leq e^{-\D^{\e}}.$$

Combining \eqref{D2l} and \eqref{D2k} with \eqref{D11b} we see that \whp,
\begin{eqnarray*}
\F_{u,i}'-\F_{u,i}&\leq& -\th(1-(k^{2k-1}+3)\e)\F_{u,i}+\binom{k-1}{i-1}\th^{k-i}\Xi_u+
\sum_{l=i+1}^{k-1}\binom{l}{i-1}\th^{l-i}\F_{u,l}+\\
&&+2k\D^{-1/2k}+2k\th^2\hp^2\brac{\binom{k-1}{i-1}\th^{k-i}\Xi_u+\sum_{l=i+1}^{k-1}\binom{l}{i-1}\th^{l-i}\F_{u,l}}\\
&\leq&\binom{k-1}{i-1}\th^{k-i}\Xi_u+\sum_{l=i+1}^{k-1}\binom{l}{i-1}\th^{l-i}\F_{u,l}-\th(1-(k^{2k-1}+3)\e)\F_{u,i}+\D^{-\e}.
\end{eqnarray*}
This confirms \eqref{B1a}.

\subsubsection{Proof of \eqref{C1}}
Fix $c$ and write
$p'=p_u'(c)=p\beta$. We consider two cases, but in both cases
$\E(\beta)=1$ and $\beta$ takes two values, 0 and $1/\Pr(\beta>0)$.
Then we have
$$\E(-p'\log p')=-p\log p -p\log(1/\Pr(\beta>0)).$$
\begin{description}
\item[(i)] $p=p_u(c)$ and $\beta=\g_u(c)/q_u(c)$ and $\g_u(c)$ is a
$\set{0,1}$ random variable with $\Pr(\beta>0)=q_u(c)$.
\item[(ii)] $p=p_u(c)=\hp$ and $\beta$ is a
$\set{0,1}$ random variable with $\Pr(\beta>0)=p_u(c)/\hp\geq
q_u(c)$.
\end{description}
Thus in both cases
$$\E(-p'\log p')\geq -p\log p-p\log 1/q_u(c).$$
Observe next that $0\leq a,b\leq 1$ implies that $(1-ab)^{-1}\leq (1-a)^{-b}$
and $-\log(1-x)\leq x+x^2$ for $0\leq x\ll 1$. So, from \eqref{2a},
\begin{eqnarray*}
\log 1/q_u(c)&\leq& -\Xi_u(c)\log(1-\th^{k-1})-\sum_{i=2}^{k-1}\F_{u,i}(c)\log(1-\th^{i-1})\\
&\leq &(\th^{k-1}+\th^{2k-2})\Xi_u(c)+\sum_{i=2}^{k-1}(\th^{i-1}+\th^{2i-2})\F_{u,i}(c).
\end{eqnarray*}
Now
\begin{eqnarray*}
\E(h_u-h_u')&\leq&-\sum_{c} p_u(c) \log p_u(c) -\E\left(\sum_{c}-p'_u(c)\log p'_u(c)\right)\\
&\le& \sum_{c} -p_u(c) \log p_u(c) -\left(\sum_{c}-p_u(c)\log p_u(c)-p_u(c)\log 1/q_u(c)\right)\\
&=& \sum_{c} p_u(c) \log 1/q_u(c)\\
&\le&(\th^{k-1}+\th^{2k-2})\sum_c p_u(c)\Xi_u(c) +\sum_c \sum_{i=2}^{k-1}(\th^{i-1}+\th^{2i-2})p_u(c)\F_{u,i}(c) \\
&=&(\th^{k-1}+\th^{2k-2})\Xi_u+\sum_{i=2}^{k-1}(\th^{i-1}+\th^{2i-2})\F_{u,i}\\
&\leq&(\th^{k-1}+\th^{2k-2})({\om}+t\D^{-\e})(1-\th/2k)^t+\sum_{i=2}^{k-1}(\th^{i-1}+\th^{2i-2})
k^{2k-2i}\om(1-\th/3k)^t\\
&\leq&k^{2k-3}\e(1-\th/3k)^t.
\end{eqnarray*}
Given the $p_u(c)$ we see that $h_u'$ is the sum of $q$ independent non-negative random variables with
values bounded by $-\hp\log\hp\leq \D^{-1/(k-1)+\e+o(1)}$. Here we have used color independence. So,
$$\Pr(h_u-h_u'\geq k^{2k-3}\e(1-\th/3k)^t+\r)\leq \exp\set{-\frac{2\r^2}{q(\hp\log\hp)^2}}=e^{-2\r^2\D^{1/k}}.$$
We take $\r=\e(1-\th/3k)^t\geq(\log\D)^{-O(1)}$ to see that
$h_u-h_u'\leq k^{2k}\e(1-\th/3k)^t$ holds \whp.
\subsubsection{Proof of \eqref{D1}}
Fix $u$ and condition on the values
$\g_w(c),\eta_w(c)$ for all $c\in C$ and all $w\notin N(u)$ and for
$w=u$. Now write $u\sim v$ to mean that $\{u,v\}$ lies in an edge of
$\ts{H}$ or some $H_i$. Then write
$$Z_u=d(u)-d'(u)\geq\frac{1}{k-1}\sum_{u\sim v}Z_{u,v}\ where \ Z_{u,v}=1_{v\notin U'}.$$
Now, for $e=uu_2\cdots u_k\in\ts{H}$ let $Z_{u,e}=\sum_{j=2}^kZ_{u,u_j}$ and if
$e=uu_2\cdots u_i\in H_i$ let $Z_{u,e}=\sum_{j=2}^iZ_{u,u_j}$.
Conditional neighborhood independence implies that the collection $Z_{u,e}$ constitute an independent set of
random variables.
Applying \eqref{hoef} to $Z_u=\sum_e Z_{u,e}$ we see that
\beq{xxx}
\Pr(Z_u\leq \E(Z_u)-\D^{2/3})\leq \exp\set{-\frac{2\D^{4/3}}{(k-1)^2\D}}=e^{-2\D^{1/3}/(k-1)^2}.
\eeq
and so we only have to estimate $\E(Z_u)$.

Fix $v\sim u$. Let $C_u(v)$ be as in \eqref{Cv1}. Condition on
$\cC$. The vertex $v$ is a member of $U'$ if none of the colors
$c\notin C_u(v))$ are tentatively activated. The activations we
consider are done independently and so
\begin{eqnarray}
\Pr(v\in U'\mid\cC)&\leq&\prod_{c\notin C_u(v)}(1-\th p_v(c)) \label{lastone}\\
&\leq&\exp\set{-\sum_{c\notin C_u(v)}\th p_v(c)}\nonumber\\
&\leq&\exp\set{-\th(1-t\D^{-\e})+\th p_v(C_u(v))}\nonumber
\end{eqnarray}
If $\cE_{\eqref{Cv1}}$ occurs then $p_v(C_u(v))\leq k^{2k}\e$.
Consequently,
$$
\Pr(v\notin U')\geq \sum_{\cC:\cE_{\eqref{Cv1}}\
occurs}\brac{1-\exp\set{-\th(1-\D^{-\e})+k^{2k}\e\th}}\Pr(\cC),$$
where the sum is well-defined due to Remark \ref{rem}. Since
$\theta\rightarrow 0$ as $\D \rightarrow \infty$, and $\e$ is
sufficiently small,
$$1-\exp\set{-\th(1-\D^{-\e})+k^{2k}\e\th}>1-e^{-\th(1-2k^{2k}\e)}>\th(1-3k^{2k}\e)>\th(1-1/k^2).$$
Recall that (\ref{Cv1}) shows that $$
 \Pr(\cE_{\eqref{Cv1}}\ fails)\le e^{-\D^{1/2k}}<1/k^2. $$
 Therefore
$$ \Pr(v \not\in U')\ge \th(1-1/k^2) \sum_{\cC:\cE_{\eqref{Cv1}}\
occurs}\Pr(\cC)>\th(1-1/k^2)^2>\th(1-1/k).$$
 This gives
$$\E(Z_u)\geq \frac1k\th d(u)$$
and \eqref{D1}.
\subsubsection{Proof of \eqref{duc}}
Observe that if $e=uu_2\cdots u_i\in H_i'\setminus H_i$ and
$\k'(e)=c$ then either

(i) there exists $1\leq j\leq k-i-1$ and vertices
$u_{i+1},\ldots,u_{i+j}$ and an edge $uu_2\cdots u_{i+j}\in H_{i+j}$
such that $u_{i+1},\ldots,u_{i+j}$ get colored in Step $t$ with $c$
and so $\g_{u_{i+1}}(c)=\cdots=\g_{u_{i+j}}(c)=1$ or

(ii) there exists $uu_2\cdots u_k\in \ts{H}$ such that
$u_{i+1},\ldots,u_k$ all receive the color $c$ and so
$\g_{u_{i+1}}(c)=\cdots=\g_{u_{k}}(c)=1$. Hence,
$$d_{i}'(u)-d_i(u)\leq \sum_{j=1}^{k-i}Z_j$$
where for $j\leq k-i-1$, $Z_j\leq Bin((1+2k\th)^t\D\hp^{k-i-j},\binom{k-i}{j}(\th\hp)^j)$ and
$Z_{k-i}\leq Bin(\D,\binom{k}{i}(\th\hp)^{k-i})$.

If $i< k-1$, then the Chernoff bound
$$Pr(Bin(n,p)\geq 2np)\leq e^{-np/3}$$
implies that for $1\leq j<k-i$,
$$\Pr\brac{Z_j\geq 2(1+2k\th)^t\binom{k-i}{j}\th^j\D\hp^{k-i}}\leq e^{-\D^\e}.$$
Similarly,
$$\Pr(Z_{k-i}\geq 2(1+2k\th)^t\th^{k-i}\D\hp^{k-i})\leq e^{-\D^\e}.$$
Therefore \whp
\begin{multline*}
d_i'(u,c)-d_i(u,c)\leq 2(1+2k\th)^t\D\hp^{k-i}\sum_{j=1}^{k-i}\binom{k-i}{j}\th^j\\
=2(1+2k\th)^t\D\hp^{k-i}((1+\th)^{k-i}-1)\leq 2k\th(1+2k\th)^t\D\hp^{k-i}.
\end{multline*}

\subsection{List Coloring}

Here we describe the small modifications needed to our argument to prove the same result for list colorings.
Each vertex $v\in V$ starts with a set $A_v$ of $2q$ available colors. Choose for each $v$ a set $B_v\subseteq A_v$
where $|B_v|=q$. Let now $C=\bigcup_{v\in V}B_v$. We initialise $p_v(c)=q^{-1}1_{c\in B_v}$ and follow the
main argument as before. When the semi-random procedure finishes, the local lemma can be used to show that
the lists $A_v\setminus B_v$ can be used to color the vertices that remain uncolored.

\end{document}